\documentclass[11pt,reqno]{amsart}
\usepackage{amsfonts,amsmath,latexsym,amssymb,bbm}    
\usepackage{graphicx,float}
\usepackage[colorlinks=true,linkcolor=blue,citecolor=blue,urlcolor=blue]{hyperref}
\usepackage{graphicx,float}
\usepackage{enumerate}
\usepackage{comment}
\usepackage{tikz}
\usepackage{longtable}
\usepackage[numbers]{natbib}
\usepackage{tabularx}

\usepackage[normalem]{ulem}

\newtheorem{notation}{Notation}[section]

\newcommand{\bu}{\mathbf{u}}

\DeclareMathOperator*{\esssup}{ess\,sup}

\numberwithin{equation}{section}

\newtheorem{theorem}{Theorem}[section]

\newtheorem{lemma}[theorem]{Lemma}
\newtheorem{rem}[theorem]{Remark}

\usepackage[top=3cm, bottom=3cm, left=2.5cm, right=2.5cm]{geometry}

\makeatletter
\@namedef{subjclassname@2020}{%
  \textup{2020} Mathematics Subject Classification}
\makeatother

\title[Pattern formation in a vasculogenesis model]{Pattern formation in a vasculogenesis model}

\author[Lahiri]{Sinchita Lahiri}
\address[S. Lahiri]{Department of Mathematics, Tulane University, New Orleans, LA 70118, USA}
\email{slahiri@tulane.edu}

\author[Zhao]{Kun Zhao}
\address[K. Zhao]{School of Mathematical Sciences, Harbin Engineering University, Harbin 150001, China}
\email{kzhao@hrbeu.edu.cn}


 
\begin{document}

\maketitle

\begin{abstract}
This paper investigates steady-state solutions of a vasculogenesis model governed by coupled partial differential equations in a bounded two-dimensional domain. Explicit steady-state solutions are analytically constructed, and their stability is rigorously analyzed under prescribed initial and boundary conditions. By employing  energy method, we prove that these solutions exhibit local asymptotic stability when specific parametric criteria are satisfied. The analysis establishes a direct connection between the stability thresholds and the system’s diffusion coefficient, offering quantitative insights into the mechanisms governing pattern formation. These results provide foundational theoretical advances for understanding self-organization in chemotaxis-driven biological systems, particularly vasculogenesis.
\end{abstract}

\section{Introduction}

\subsection{Background} 

Vasculogenesis, distinct from angiogenesis (the sprouting of new vessels from pre-existing ones \cite{3}), refers to the {\it de novo} formation of blood vessels during embryonic development \cite{9}. This process involves endothelial progenitor cells (angioblasts) that differentiate, migrate, and self-assemble into primitive vascular networks. Biologically, vasculogenesis progresses through four distinct phases:
\begin{itemize}
\item[1.] Chemotactic migration: Cells migrate toward regions of higher cellular density until colliding with neighboring cells;
\item[2.] Topological networking: Post-collision adhesion establishes continuous multicellular connections;
\item[3.] Mechanistic remodeling: The network undergoes collective motion and substrate-dependent thinning;
\item[4.] Morphological maturation: Individual cells fold to form capillary lumens, completing vascular network architecture.
\end{itemize}
A comprehensive understanding of these mechanisms holds significant implications for developmental biology, regenerative medicine, and disease therapeutics.

In mathematical biology, vasculogenesis is commonly modeled through partial differential equations (PDEs) that integrate cell motility, chemical signaling, and substrate interactions. This study focuses on the following nonlinear PDE system proposed by Gamba et al.\,\cite{8} to replicate {\it in vitro} vasculogenesis dynamics:
\begin{subequations}\label{eqn:vasculogenesismodel}
\begin{alignat}{3}
\rho_t + \nabla \cdot(\rho \mathbf{u}) &=0, \label{vas1}\\
(\rho \mathbf{u})_t + \nabla \cdot(\rho \mathbf{u}\otimes \mathbf{u})+\nabla P(\rho) &=-\alpha\rho \mathbf{u}+\beta\rho\nabla\phi, \label{vas2}\\
\tau\phi_t - d\Delta\phi + a\phi &= b\rho, \label{vas3}
\end{alignat}
\end{subequations}
where $\mathbf{x}\in \mathbb{R}^N$ and $t>0$. For future reference, we refer to this model as the Fluid-Reaction-Diffusion system, or simply the FRD system. The model employs a phenomenological framework grounded in six key assumptions:
\begin{itemize}
\item[1.] Cell population conservation (no proliferation/apoptosis, per the continuity equation for $\rho$);
\item[2.] Motion persistence (convective term $\nabla \cdot(\rho \mathbf{u}\otimes \mathbf{u})$);
\item[3.] Inviscid Newtonian fluid behavior (pressure term $\nabla P(\rho)$);
\item[4.] Frictional damping (term $-\alpha\rho\mathbf{u}$);
\item[5.] Chemotactic coupling (nonlocal term $\beta\rho\nabla\phi$ for growth factor $\phi$);
\item[6.] Growth factor dynamics (reaction-diffusion equation with secretion rate $b>0$, degradation $a>0$, diffusion $d>0$, and relaxation time $\tau\geqslant 0$).
\end{itemize}
Here, the unknown functions $\rho$, $\mathbf{u}$, and $\phi$ denote the density of cell population, velocity of ensemble of cells, and concentration of growth factor at position $\mathbf{x}$ at time $t$, respectively. The parameter $\alpha>0$ quantifies damping effects, whereas $\beta$ determines chemotactic polarization: $\beta>0$ represents attraction, and $\beta<0$ corresponds to repulsion. Gamba et al.\,\cite{8} demonstrated this model's capacity to recapitulate the first three vasculogenic stages (chemotaxis, network formation, remodeling), with extended theoretical foundations detailed by Ambrosi et al.\,\cite{1}.

\subsection{Literature Review}

Immediately following the design of the FRD system \eqref{eqn:vasculogenesismodel}, rigorous mathematical analyses were conducted to investigate the model's qualitative properties, including global well-posedness, long-time solution behavior, connections to Keller-Segel-type chemotaxis models, and pattern solution stability. Below we systematically review relevant studies categorized by spatial domain geometry (with/without physical boundaries).

\underline{Cauchy problem in whole space}: The first mathematical result was established in \cite{12} for the ``viscous" model through introduction of the $\Delta\mathbf{u}$ term into \eqref{vas2}. This work demonstrated linear stability of constant ground states $(\bar{\rho},0,\bar\phi)$ ($\bar\rho,\bar\phi>0$) under the condition:
\begin{equation}\label{eqn:pressure condition}
  bP'(\bar{\rho}) - a\alpha\bar{\rho}>0.
\end{equation}
Subsequent asymptotic analysis by \cite{5} revealed convergence of \eqref{eqn:vasculogenesismodel} to Keller-Segel model as damping coefficient $\alpha\to\infty$, later mathematically substantiated in \cite{7}. Numerical comparisons between model \eqref{eqn:vasculogenesismodel} and Keller-Segel systems were further explored in \cite{16}. For small perturbations of constant ground states ($\|(\rho_0-\bar\rho,\mathbf{u}_0,\phi_0-\bar\phi)\|_{H^s(\mathbb{R}^N)}\ll 1$, $s>N/2+1$) with sufficiently small $\bar{\rho}>0$,  global well-posedness of vacuum-free classical solutions was proven in \cite{19,20}, with algebraic decay $(1+t)^{-3/4}$ toward equilibrium. Notably, \cite{14} successfully eliminated the small-density constraint under \eqref{eqn:pressure condition}, establishing the global well-posedness of classical solutions with enhanced decay rates in $\mathbb{R}$. Their subsequent work \cite{15} demonstrated local asymptotic stability of diffusion waves in $\mathbb{R}^3$. Recent progress by \cite{6} extended well-posedness results to hybrid Besov spaces. More recently, the existence/non-existence of stationary vacuum solutions (bump-type profiles) in $\mathbb{R}^2$ was established in the radially symmetric setting \cite{HHYZ}.

\underline{Boundary effects}: The incorporation of physical boundaries substantially complicates the mathematical analysis of the FRD model, resulting in comparatively limited mathematical investigations relative to the whole-space Cauchy problem. Notable exceptions include construction of stationary vacuum solutions (bump-type profiles) in bounded intervals under zero-flux boundary conditions \cite{2,4}, and local stability of transition layers on the semi-infinite domain $\mathbb{R}^+$ \cite{10}.

\subsection{Open question and current work}

Recent work by \cite{18} established the first rigorous mathematical results on the existence and stability of non-constant stationary solutions to \eqref{eqn:vasculogenesismodel} in \emph{bounded domains}. However, these findings are currently restricted to one-dimensional space, despite the fact that multi-dimensional settings are far more biologically relevant. Indeed, {\it in vivo}, {\it in vitro}, and {\it in silico} studies of vasculogenesis inherently occur in higher-dimensional spaces. To bridge this gap, we extend the analysis to two-dimensional bounded domains, constructing explicit steady states of system \eqref{eqn:vasculogenesismodel} and investigating their long-time asymptotic stability. Specifically, we examine system \eqref{eqn:vasculogenesismodel} under the functional setting  $P(\rho)=A_0\rho^\gamma$ (where $A_0>0$ and $\gamma\geqslant 1$ are constants) on the unit square domain $\Omega=(0,1)\times(0,1)$, subject to the following initial and boundary conditions:
\begin{subequations}\label{eqn:ourmodel}
\begin{alignat}{2} 
\bu(\mathbf{x},t)\cdot\mathbf{n}&=0, \quad & \mathbf{x}&\in \partial\Omega,\ t\geqslant 0, \label{bc1}\\
\phi(\mathbf{x},t)&=0, \quad & \mathbf{x}&\in \partial\Omega,\ t\geqslant 0,\quad \text{or} \label{bc2}\\
\nabla\phi(\mathbf{x},t)\cdot \mathbf{n}&=0,\quad & \mathbf{x}&\in \partial\Omega,\ t\geqslant 0, \label{bc3}
\end{alignat}
\end{subequations}
where $\mathbf{n}$ denotes the unit outward normal to $\partial\Omega$.

We remark that on one side, when $\gamma=2$, the quadratic pressure function  aligns with previous studies on steady-state solutions of the FRD model \cite{2,4,18}. This specific functional form leads to linear steady-state equations, facilitating the derivation of explicit solutions, a crucial advantage for analytical tractability. Additionally, the square domain’s geometry is instrumental in constructing these solutions, as it permits the use of eigenfunction expansion method.
On the other side, our analysis focuses on the two-dimensional setting for two key reasons:
\begin{itemize}
\item[1.] Model validation: Current three-dimensional simulations of the model cannot yet replicate essential features of {\it in vivo} or {\it in vitro} vasculogenesis.
\item[2.] Theoretical constraint: The validity of foundational analytical tools, such as Sobolev embeddings, Gagliardo-Nirenberg inequalities, and elliptic regularity estimates, remains unverified for three-dimensional polyhedral domains (e.g., cuboids).
\end{itemize}
Furthermore, the no-normal-flow (non-characteristic) boundary condition for the velocity field is standard for inviscid fluid flows, while the Dirichlet and Neumann boundary conditions for $\phi$ follow the same formulations as in \cite{18} and \cite{2,4}, respectively.

The rest of this paper is organized as follows. In Section 2, we construct and analyze the stability of steady states for system \eqref{eqn:vasculogenesismodel} under quadratic pressure, $P(\rho)=A_0\rho^2$, subject to boundary conditions \eqref{bc1} $\&$ \eqref{bc2}. In Section 3, we extend our analysis to the case of ``$\gamma$-law pressure'', $P(\rho)=A_0\rho^\gamma$ ($\gamma\geqslant 1$), with the boundary conditions \eqref{bc1} $\&$ \eqref{bc3}. The paper ends with concluding remarks and future directions in Section 4.

\section{Steady state and its stability under \eqref{bc1} $\&$ \eqref{bc2}}

In this section, we investigate the steady-state solution for \eqref{eqn:vasculogenesismodel} with quadratic pressure:
\begin{subequations}\label{HP-A}
\begin{alignat}{3}
\rho_t + \nabla \cdot(\rho \mathbf{u}) &=0, \label{HP-A1}\\
(\rho \mathbf{u})_t + \nabla \cdot(\rho \mathbf{u}\otimes \mathbf{u})+\nabla (A_0\rho^2) &=-\alpha\rho \mathbf{u}+\beta\rho\nabla\phi, \label{HP_a2}\\
\tau\phi_t - d\Delta\phi + a\phi &=b\rho, \label{HP-A3}
\end{alignat}
\end{subequations}
under the boundary conditions \eqref{bc1} and \eqref{bc2}.

\subsection{Finding steady state} 

Due to the dissipation mechanism induced by linear damping and the boundary condition imposed on $\bu$, it is natural to anticipate that the equilibrium velocity is zero. For the other components of the steady state, let us denote them by $\hat{\rho}$ and $\hat{\phi}$. Then, these functions satisfy the following boundary value problem:
\begin{subequations}\label{eqn:steady_state}
\begin{alignat}{3}
2A_0\hat{\rho}\nabla\hat{\rho} - \beta\hat{\rho}\nabla\hat{\phi} &=0, \quad &\mathbf{x}& \in\Omega, \label{SS1}\\
d\Delta\hat{\phi}-a\hat{\phi}+b\hat{\rho} &=0, &\mathbf{x}&\in\Omega; \label{SS2}\\
\hat{\phi}&=0, &\mathbf{x}&\in\partial\Omega. \label{SS3}
\end{alignat}
\end{subequations}
Note that because of cellular mass conservation and the boundary condition imposed on $\phi$, the equilibrium chemical concentration $\hat\phi$ cannot be a constant. Meanwhile, our primary concern is the existence and stability of vacuum-free steady state. Under this consideration, \eqref{SS1} leads to 
\begin{align}\label{SS4}
\hat{\rho}(\mathbf{x})=\frac{\beta}{2A_0}\hat{\phi}(\mathbf{x})+\frac{\hat{C}}{2A_0},
\end{align}
where the constant $\hat{C}$ is to be determined later. Substituting \eqref{SS4} into \eqref{SS2}, we obtain
\begin{equation}\label{eqn:solvephi}
\Delta\hat{\phi}(\mathbf{x}) - \Lambda \hat\phi(\mathbf{x}) +\hat{D}=0,
\end{equation}
where the constants are defined by
\begin{align}\label{SS5}
\Lambda \triangleq \frac{2aA_0-b\beta}{2dA_0}\quad \text{and}\quad \hat{D} \triangleq \frac{b\hat{C}}{2dA_0}.
\end{align} 
Applying the method of eigenfunction expansion, we obtain the solution to \eqref{eqn:solvephi} subject to the boundary condition \eqref{SS3} in the form:\begin{align}\label{SS7}
\hat{\phi}(\mathbf{x})=\frac{8b\hat{C}}{dA_0 \pi^2}\sum_m\sum_n \frac{\sin(m\pi x)\sin(n\pi y)}{mn[(m^2+n^2)\pi^2+\Lambda]}, \quad m,n \in\mathbb{N}\ \text{are odd},
\end{align}  
where $\mathbf{x}=(x,y)\in\Omega$. To ensure the validity of the solution, we must require that for any odd natural numbers $m$ and $n$, $\Lambda\neq -(m^2+n^2)\pi^2$. Substituting \eqref{SS7} into \eqref{SS4} yields
\begin{equation}\label{SS8}
\hat{\rho}(\mathbf{x})=\frac{4b\beta\hat{C}}{dA_0^2\pi^2}\sum_m\sum_n\frac{\sin(m\pi x)\sin(n\pi y)}{mn[(m^2+n^2)\pi^2+\Lambda]} +\frac{\hat{C}}{2A_0}, \quad \mathbf{x} \in\Omega.
\end{equation}
To determine the constant $\hat{C}$, we invoke the conservation of total cellular mass, which holds under the no-normal-flow boundary condition. This yields the constraint:
\begin{align}\label{SS9}
\int_\Omega \hat{\rho}(\mathbf{x})\mathrm{d}\mathbf{x} = \frac{4b\beta\hat{C}}{dA_0^2\pi^4}\sum_m\sum_n\frac{1}{m^2n^2[(m^2+n^2)\pi^2+\Lambda]}+\frac{\hat{C}}{2A_0} = \int_\Omega \rho_0(\mathbf{x}) \mathrm{d}\mathbf{x},
\end{align}
from which we deduce that 
\begin{align}\label{SS10}
\hat{C}=
\left(\int_\Omega \rho_0(\mathbf{x}) \mathrm{d}\mathbf{x}\right) \left[\frac{4b\beta}{dA_0^2\pi^4} \sum_m\sum_n\frac{1}{m^2n^2[(m^2+n^2)\pi^2+\Lambda]}+\frac{1}{2A_0}\right]^{-1}.
\end{align}
Note that for this formula to hold, the bracketed expression must be non-zero.


\subsection{Stability of steady state} 

We now prove the long-time asymptotic stability of the steady-state solution constructed above. At present, the coupled hyperbolic-parabolic structure of the general system \eqref{HP-A}, subject to boundary conditions \eqref{bc1} and \eqref{bc2}, introduces significant analytical challenges that exceed the scope of our current methodology. The primary difficulty arises from the boundary condition imposed on $\phi$: under \eqref{bc1} and \eqref{bc2}, no boundary condition can be explicitly derived for $\rho$, in contrast to the purely hydrodynamic case. This absence is critical because energy estimates rely on integration by parts, which requires boundary terms for $\rho$ to control spatial derivatives of the solution. Consequently, we are compelled to restrict our analysis to temporal derivatives. For the hyperbolic-parabolic system \eqref{HP-A}, given initial data $\rho_0,\bu_0 \in H^k(\Omega)$, the $H^k$-estimate of $(\rho,\bu)$ depends on $\|(\rho_0,\bu_0)\|_{H^{2k-2}}$. This dependency confines the energy estimates to the case $k=2$. However, establishing well-posedness for classical solutions in dimension two necessitates an energy space $H^k(\Omega)$ with $k\geqslant 3$. Thus, our approach cannot be utilized to prove the stability of steady-state solutions for \eqref{HP-A}. Due to such a technical limitation, our analysis is restricted to the fast-relaxation case ($\tau=0$) of system \eqref{HP-A}:
 \begin{subequations}\label{AA1}
\begin{alignat}{3}
\rho_t + \nabla \cdot(\rho \mathbf{u}) &=0, \label{aa1}\\
(\rho \mathbf{u})_t + \nabla \cdot(\rho \mathbf{u}\otimes \mathbf{u})+ \nabla (A_0\rho^2) &=-\alpha\rho \mathbf{u}+\beta\rho\nabla\phi, \label{aa2}\\
-d\Delta\phi + a\phi&=b\rho.\label{aa3}
\end{alignat}
\end{subequations}
The stability analysis for system \eqref{eqn:vasculogenesismodel} will be addressed in Section 3, where we consider Neumann boundary condition for $\phi$. For now, we focus on presenting the stability result for system \eqref{AA1}. To this end, we first introduce the following notation for clarity and future reference.

\begin{notation} 
Throughout this work, we employ the following notation: $\|\cdot\| \triangleq \|\cdot\|_{L^2(\Omega)}$ denotes the standard $L^2$ norm; $\|\cdot\|_{H^s}$ and $\|\cdot\|_{L^\infty}$ represent the norms of $H^s(\Omega)$ and $L^{\infty}(\Omega)$, respectively. For function $f$ and $s\in\mathbb{N}$, we define the total energy of order $s$ as: $\|f(t)\|_s^2 \triangleq \sum_{k=0}^s \|(\partial_t^kf)(t)\|_{H^{s-k}}^2$. For vector-valued functions $(f_1,...,f_n)$, we use the convention: $\|(f_1,...,f_n)\|^2_* \triangleq\sum_{i=1}^n \|f_i\|_*^2$, where $*$ may denote $L^2$, $H^s$, $L^{\infty}$, or $s$. Unless otherwise stated, $C$ denotes a generic positive constant independent of the solution functions, whose value may change from line to line.
\end{notation}

\begin{theorem}\label{main}
Consider the initial-boundary value problem of \eqref{AA1} with initial data $(\rho_0,\bu_0)$ and boundary conditions \eqref{bc1} and \eqref{bc2}. Let $A_0,\alpha, \beta, d, a, b>0$ be constants, and let $(\hat\rho,\hat\phi)$ be given by \eqref{SS7}, \eqref{SS8}, and  \eqref{SS10}. Suppose that there is a large constant $d_0>0$ such that $d \geqslant d_0$, and that the initial data satisfy the following conditions:
\begin{itemize}
\item $\rho_0>0$ and $(\rho_0,\bu_0) \in [H^3(\Omega)]^3$;

\item $(\rho_0,\bu_0)$ are compatible with the boundary conditions;

\item there is a small constant $\delta_0>0$ such that $\|(\rho_0-\hat{\rho},\bu_0)\|^2_{H^3}\leqslant \delta_0$.

\end{itemize}
Then there exists a unique solution to the initial-boundary value problem, such that for $\forall t>0$,
\begin{equation*}
\begin{aligned}
 \mathbb{E}(t) + \int_0^t \mathbb{E}(\tau)\mathrm{d}\tau \leqslant C,
\end{aligned}
\end{equation*}
where $\mathbb{E}(t)\triangleq \|(\tilde{\rho}, \bu)(t)\|^2_3 + \sum_{k=0}^3 \|(\partial _t^k \tilde{\phi})(t)\|^2_{H^{5-k}}$, $\tilde{\rho}=\rho-\hat{\rho}$, and $\tilde{\phi}=\phi-\hat{\phi}$. Moreover, there are positive constants $\eta_1$ and $\eta_2$, such that for $\forall t>0$, $\mathbb{E}(t) \leqslant \eta_1 e^{-\eta_2 t}$.
\end{theorem}

\begin{rem}
We emphasize that a sufficiently large diffusion coefficient $d$ is crucial for ensuring both the biological relevance of the steady state and the validity of its stability analysis. First, when $d$ is large, the modulus of the constant $\Lambda$
defined in \eqref{SS5} becomes small, which guarantees the existence and uniqueness of the steady state given by \eqref{SS7} and \eqref{SS8}. This observation has an important biological implication: in systems with chemo-attraction ($\beta>0$) and chemical production, the chemotactic sensitivity of cells must remain bounded. Such a constraint is biologically reasonable, as an excessively strong chemotactic response, particularly under sustained chemical production, could lead to uncontrolled cellular aggregation. This would disrupt the formation of structured vascular networks, where precise spatial organization and regulated migration are essential. Second, a large diffusion coefficient ensures that the constant $\hat{C}$ in \eqref{SS10} is strictly positive, thereby preserving the positivity and boundedness of the equilibrium density. Combined with \eqref{SS2}, \eqref{SS3}, and the maximum principle, this guarantees the positivity of $\hat\phi$, reinforcing the biological plausibility of the steady-state solution. Furthermore, our asymptotic analysis will demonstrate that maintaining $\hat\rho$ strictly positive and uniformly bounded is critical for the stability of the steady state. Due to the complexity of the expression for $d_0$, we defer its detailed derivation to the proof.
\end{rem}

As usual, the proof of Theorem \ref{main} is structured into three integral parts: the establishment of local well-posedness, the derivation of {\it a priori} estimates, and the application of a continuation argument. The initial segment, concerning local well-posedness, is substantiated by classic reasoning presented in \cite{21}, complemented by the energy estimates developed in this paper. For the sake of conciseness, we present the outcome of this analysis, omitting the intricate technicalities.

\begin{lemma}\label{local}
Under the conditions of  Theorem \ref{main}, there exists some finite $T\in (0,\infty)$ and a unique solution to the initial-boundary value problem of \eqref{AA1} subject to the boundary conditions \eqref{eqn:ourmodel}, such that $(\rho,\bu) \in [L^\infty(0,T;H^3(\Omega))]^3$ and $\phi \in L^\infty(0,T;H^5(\Omega))$. 
\end{lemma}

The core of this section is devoted to establishing the {\it a priori} estimates for the local solution derived in Lemma \ref{local}. Once these {\it a priori} estimates are secured, the global well-posedness and long-time behavior of the solution can be deduced through a standard continuation argument. 
Our proof of the {\it a priori} estimates proceeds via the following key steps:
\begin{itemize}
\item[1.] Reformulation via a sound-speed transformation.
\item[2.] Reduction of spatial derivatives to temporal derivatives.
\item[3.] Energy estimates for temporal derivatives.
\item[4.] Recovery of density dissipation.
\item[5.] Coupling of hyperbolic and elliptic estimates.
\end{itemize}
While the first four steps have appeared in prior studies of damped Euler equations \cite{17}, the present analysis introduces significant additional complexity due to the nonlocal coupling between chemical concentration and fluid variables. In addition, as previously noted, the sufficiently large diffusion coefficient proves essential in the final step when closing the full energy scheme. Furthermore, in the process of deriving the {\it a priori} estimates, we extensively employ fundamental analytical tools, including the Poincar\'{e} inequality, Sobolev embedding theorems, and elliptic regularity estimates. The applicability of these tools in planar polygonal domains, particularly in the unit square, is rigorously justified by a body of prior works \cite{Bar,Gri,Lac}.

\vspace{.1 in}

{\bf Preparation}:
To prove Theorem \ref{main}, we first reformulate the system of equations,  \eqref{AA1}, via the sound-speed transformation: $\sigma \triangleq 2\sqrt{2A_0\rho}$. In terms of the transformed density $\hat\sigma \triangleq 2 \sqrt{2 A_0 \hat{\rho}}$ and the perturbed variables $\tilde{\sigma} \triangleq \sigma - \hat{\sigma}$ and $\tilde{\phi} \triangleq \phi -\hat{\phi}$, system \eqref{AA1} becomes
\begin{subequations}\label{AA2}
\begin{alignat}{3}
2 \tilde{\sigma}_t + 2 \bu\cdot \nabla{\tilde{\sigma}} + \tilde{\sigma} \nabla\cdot \bu + 2 \bu\cdot \nabla\hat{\sigma} + \hat{\sigma} \nabla\cdot\bu &=0, \label{aa4}\\
2 \bu_t + 2 (\bu\cdot \nabla )\bu + \tilde{\sigma}\nabla{\tilde{\sigma}} + \tilde{\sigma} \nabla{\hat{\sigma}} + \hat{\sigma}\nabla{\tilde{\sigma}}  
&=-2 \alpha \bu + 2 \beta \nabla{\tilde{\phi}}, \label{aa5}\\
8 d A_0 \Delta\tilde{\phi}- 8a A_0 \tilde{\phi} + b (\tilde{\sigma} + 2 \hat{\sigma})\tilde{\sigma} &= 0. \label{aa6}
\end{alignat}
\end{subequations}
Additionally, the initial and boundary conditions turn into 
\begin{subequations}\label{AA3}
\begin{alignat}{2}
(\tilde{\sigma}, \bu)(\mathbf{x},0) &= (2\sqrt{2 A_0 \rho_0}- 2 \sqrt{2 A_0 \hat{\rho}}, \bu_0)(\mathbf{x}), \label{aa7}\\
\bu\cdot\mathbf{n}|_{\partial \Omega}&=0=\tilde{\phi}|_{\partial \Omega}. \label{aa8}
\end{alignat}
\end{subequations}
This completes the first step outlined above concerning the proof of Theorem \ref{main}. In the subsequent subsections, we shall frequently utilize the {\it a priori} assumption:
\begin{equation}\label{AA4}
\esssup_{t \in [0,T]} \mathtt{E}(t)\triangleq \esssup_{t \in [0,T]} \|(\tilde{\sigma},\bu)(t)\|^2_3 \leqslant \delta,
\end{equation} 
where $T > 0$ denotes the lifespan of the local solution and $\delta>0$  is a small constant. The smallness of $\delta$ can be realized through the smallness assumption of the initial perturbation in Theorem \ref{main} and the local well-posedness theory.

\vspace{.1 in}

{\bf Reduction of total energy:} To reduce the total energy, we define the following quantities:
\begin{align}\label{AA5}
\mathtt{U}_1(t) \triangleq \sum_{k=1}^3\|(\partial_t^k \tilde{\sigma})(t)\|^2+\sum_{k=0}^3 \|(\partial_t^k \bu)(t)\|^2 \quad \text{and} \quad \mathtt{V}_1(t)\triangleq  \|\omega(t)\|_{2}^2.
 \end{align}
 where $\omega=\nabla\times\bu$ denotes the 2D vorticity. 
 
{\bf Step 1.} The initial step of the reduction process is to establish a Poincar\'e-type inequality for $\tilde\sigma$, i.e., $\|\tilde\sigma\| \leqslant C\|\nabla\tilde\sigma\|$. This follows from the same arguments as in \cite{18}, by taking advantage of the Poincar\'e inequality for $\rho-\hat\rho$ (as it is mean-free), the boundedness of $\hat\rho$ given by \eqref{SS14}, and the {\it a priori} assumption \eqref{AA4}. We omit the technical details for brevity.

\textbf{Step 2.} Taking the $L^2$ inner product of \eqref{aa6} with $-\tilde{\phi}$, we can show that
\begin{align}\label{AA6}
8dA_0\|\nabla\tilde\phi\|^2+8aA_0\|\tilde\phi\|^2=b\int_\Omega (\tilde\sigma+2\hat\sigma)\tilde\sigma \tilde\phi\mathrm{d}\mathbf{x} &\leqslant b\|(\tilde{\sigma} + 2 \hat{\sigma})\|_{L^{\infty}}\|\tilde{\sigma}\|\|\tilde{\phi}\|\notag\\
&\leqslant C\|\tilde{\sigma}\|\|\tilde{\phi}\|,
\end{align}
where the {\it a priori} assumption \eqref{AA5} and boundedness of $\hat\sigma$ are applied when deriving the last inequality. By Young's inequalities, we can improve the last term in \eqref{AA6} as 
\begin{align}\label{AA7}
C\|\tilde{\sigma}\|\|\tilde{\phi}\| \leqslant 4aA_0\|\tilde\phi\|^2 + C\|\tilde\sigma\|^2.
\end{align}
Substituting \eqref{AA7} into \eqref{AA6} yields
\begin{align}\label{AA8}
8dA_0\|\nabla\tilde\phi\|^2+4aA_0\|\tilde\phi\|^2 \leqslant C\|\tilde\sigma\|^2,
\end{align}
which implies 
\begin{align}\label{AA9a}
\|\nabla\tilde\phi\|^2 \leqslant \frac{C}{d}\|\tilde\sigma\|^2 \leqslant \frac{C}{d}\|\nabla\tilde\sigma\|^2,
\end{align}
where we applied the Poincar\'e-type inequality for $\tilde\sigma$. By elliptic regularity estimates, we have 
\begin{align}\label{AA9b}
\|\tilde\phi\|_{H^3}^2 \leqslant \frac{C}{d^2}\big(\|\tilde\phi\|_{H^1}^2 + \|(\tilde\sigma+2\hat\sigma)\tilde\sigma\|_{H^1}^2\big).
\end{align}
Utilizing Sobolev embedding and Poincar\'e's inequality, we can show that 
\begin{align}\label{AA9c}
\|(\tilde\sigma+2\hat\sigma)\tilde\sigma\|_{H^1}^2 \leqslant C\big(\|\nabla\tilde\sigma\|_{H^2}^2 + \|\hat\sigma\|_{H^3}^2\big) \|\nabla\tilde\sigma\|^2.
\end{align}
Combining \eqref{AA9a}, \eqref{AA9b}, and \eqref{AA9c}, we obtain
\begin{align}\label{AA9}
\|\tilde\phi\|_{H^3}^2 \leqslant \frac{C}{d^2}\big(1 + \|\nabla\tilde\sigma\|_{H^2}^2 + \|\hat\sigma\|_{H^3}^2\big) \|\nabla\tilde\sigma\|^2.
\end{align}
We emphasize that \eqref{AA9} plays a critical role in reducing the influence of nonlocal coupling between the chemical concentration and fluid variables.

{\bf Step 3.} From \eqref{aa4} and \eqref{aa5}, we have
\begin{align}
\nabla\cdot\bu &=-\frac{1}{\tilde{\sigma}+\hat{\sigma}}\left (2 \tilde{\sigma}_t+2 \bu \cdot\nabla{\tilde{\sigma}}+2 \bu\cdot\nabla\hat{\sigma} \right ), \label{AA10}\\
  \nabla{\tilde{\sigma}} &=-\frac{1}{\tilde{\sigma}+\hat{\sigma}}\big(2 \bu_t+2 (\bu\cdot \nabla )\bu +\tilde{\sigma} \nabla{\hat{\sigma}} + 2 \alpha \bu - 2 \beta \nabla{\tilde{\phi}}\big). \label{AA11}
\end{align}
Upon taking $L^2$ inner product and applying \eqref{AA4} and boundedness of $\hat\sigma$, we can show that
\begin{align}
\|\nabla\cdot\bu\|^2 &\leqslant C\big(\|\tilde\sigma_t\|^2+\|\bu\|_{L^\infty}^2\|\nabla\tilde\sigma\|^2+ \|\bu\|^2\|\nabla\hat\sigma\|_{L^\infty}^2\big), \label{AA12}\\
\|\nabla\tilde\sigma\|^2 &\leqslant C\big(\|(\bu_t,\bu)\|^2+\|\bu\|_{L^\infty}^2\|\nabla\bu\|^2+ \|\nabla\tilde\sigma\|^2\|\nabla\hat\sigma\|_{L^\infty}^2+\|\nabla\tilde\phi\|^2\big), \label{AA13}
\end{align}
where the Poincar\'e-type inequality for $\tilde\sigma$ is applied. Now using the div-curl estimate (cf.\,\cite{22}):
\begin{align}\label{AA14}
 \|\bu\|_{H^s}\leqslant C\left(\|\nabla\times\bu\|_{H^{s-1}}+\|\nabla \cdot\bu\|_{H^{s-1}} + \|\bu\|_{H^{s-1}}\right),
\end{align}
as well as the definition of $\mathtt{E}(t)$ given by \eqref{AA4}, we update \eqref{AA12} as 
\begin{align}
\|\bu\|_{H^1}^2 \leqslant C\big(\|(\tilde\sigma_t,\bu)\|^2 + \|\omega\|^2 + \mathtt{E}^2+ \|\bu\|^2\|\nabla\hat\sigma\|_{L^\infty}^2\big). \label{AA15}
\end{align}
Meanwhile, inserting \eqref{AA9a} into \eqref{AA13} and applying the Poincar\'e-type inequality for $\tilde\sigma$ yield
\begin{align}
\|\tilde\sigma\|_{H^1}^2 \leqslant C\big(\|(\bu_t,\bu)\|^2+\mathtt{E}^2+ \|\nabla\tilde\sigma\|^2\|\nabla\hat\sigma\|_{L^\infty}^2+ d^{-1}\|\nabla\tilde\sigma\|^2\big). \label{AA16}
\end{align}
Since $\nabla\hat\sigma \to 0$ as $d\to\infty$, when $d$ is sufficiently large, we further upgrade \eqref{AA15} and \eqref{AA16} as
 \begin{align}
\|\bu\|_{H^1}^2 &\leqslant C\big(\|(\tilde\sigma_t,\bu)\|^2 + \|\omega\|^2+ \mathtt{E}^2\big), \label{AA17}\\
\|\tilde\sigma\|_{H^1}^2 &\leqslant C\big(\|(\bu_t,\bu)\|^2+\mathtt{E}^2\big). \label{AA18}
\end{align}
In a similar fashion, by taking temporal derivatives of \eqref{AA10}--\eqref{AA11}, we can show that 
 \begin{align}
\|\partial_t^k\bu\|_{H^1}^2 &\leqslant C\big(\|\partial_t^k(\tilde\sigma_t,\bu)\|^2 + \|\partial_t^k\omega\|^2+ \mathtt{E}^2\big), \quad k=1,2,\label{AA19}\\
\|\partial_t^k\tilde\sigma\|_{H^1}^2 &\leqslant C\big(\|\partial_t^k(\bu_t,\bu)\|^2+\mathtt{E}^2\big), \qquad k=1,2.\label{AA20}
\end{align}
where we used similar estimates as \eqref{AA9a} for $\tilde\phi_t$ and $\tilde\phi_{tt}$, whose derivation is omitted for brevity.

{\bf Step 4.} By taking spatial derivatives of \eqref{AA10}--\eqref{AA11}, it can be shown that
 \begin{align}
\|\bu\|_{H^2}^2 &\leqslant C\big(\|(\tilde\sigma_t,\bu)\|_{H^1}^2 + \|\omega\|_{H^1}^2 + \mathtt{E}^2\big) \leqslant C\big(\|(\tilde\sigma_t,\bu_{tt},\bu_t,\bu)\|^2 + \|\omega\|_{H^1}^2+ \mathtt{E}^2\big), \label{AA21}\\
\|\tilde\sigma\|_{H^2}^2 &\leqslant C\big(\|(\bu_t,\bu)\|_{H^1}^2+\mathtt{E}^2\big) \leqslant C\big(\|(\tilde\sigma_{tt},\tilde\sigma_t,\bu_t,\bu)\|^2+\|(\omega_t,\omega)\|^2+\mathtt{E}^2\big), \label{AA22}
\end{align}
where \eqref{AA17}--\eqref{AA20} are applied. A further iteration on the temporal derivative yields 
 \begin{align}
\|\partial_t\bu\|_{H^2}^2 &\leqslant C\big(\|(\tilde\sigma_{tt},\bu_{ttt},\bu_{tt},\bu_t)\|^2 + \|\omega_t\|_{H^1}^2+ \mathtt{E}^2\big), \label{AA23}\\
\|\partial_t\tilde\sigma\|_{H^2}^2 &\leqslant C\big(\|(\tilde\sigma_{ttt},\tilde\sigma_{tt},\bu_{tt},\bu_t)\|^2+\|(\omega_{tt},\omega_t)\|^2+\mathtt{E}^2\big). \label{AA24}
\end{align}
Finally, taking the second order spatial derivatives of  \eqref{AA10}--\eqref{AA11}, we can derive the following:
 \begin{align}
\|\bu\|_{H^3}^2 &\leqslant C\big(\|(\tilde\sigma_{ttt},\tilde\sigma_{tt},\tilde\sigma_t,\bu_{tt},\bu_{t},\bu)\|^2 + \|(\omega_{tt},\omega_t)\|^2+ \|\omega\|_{H^2}^2+ \mathtt{E}^2\big), \label{AA25}\\
\|\tilde\sigma\|_{H^3}^2 &\leqslant C\big(\|(\tilde\sigma_{tt},\tilde\sigma_t,\bu_{ttt},\bu_{tt},\bu_t,\bu)\|^2+\|(\omega_t,\omega)\|_{H^1}^2+\mathtt{E}^2\big). \label{AA26}
\end{align}
The combination of \eqref{AA19}--\eqref{AA20} and \eqref{AA23}--\eqref{AA26} leads to 
 \begin{align}
\|(\bu,\tilde\sigma)\|_{3}^2 \leqslant C\Big(\sum_{k=1}^3\|\partial_t^k\tilde\sigma\|^2+\sum_{k=0}^3\|\partial_t^k\bu\|^2 + \|\omega\|_2^2 + \mathtt{E}^2\Big). \label{AA27}
\end{align}
 In view of \eqref{AA4}, \eqref{AA5}, and \eqref{AA27}, we see that 
 \begin{align}
\mathtt{E} \leqslant C\big(\mathtt{U}_1+\mathtt{V}_1+\mathtt{E}^2\big).\label{AA28}
\end{align}
When $\tilde\delta$ is sufficiently small, we upgrade \eqref{AA28} as 
 \begin{align}
\mathtt{E} \leqslant C\big(\mathtt{U}_1+\mathtt{V}_1\big).\label{AA29}
\end{align}
This completes the reduction of total energy to temporal derivatives.

\vspace{.1 in}

{\bf Energy estimates of temporal derivatives}: The estimates of the temporal derivatives of the solution are obtained in a standard manner by exploiting the symmetric structure of \eqref{aa4}--\eqref{aa5}, the Sobolev embedding theorems, the boundedness of the steady-state solution, and elliptic regularity estimates. For brevity, we provide only the detailed derivation for the zeroth-order estimate, leaving the higher-order cases to interested readers.

{\bf Step 1.} Taking $L^2$ inner product of \eqref{aa4} with $\tilde{\sigma}$, we have
\begin{align}\label{AA30}
\frac{\mathrm{d}}{\mathrm{d}t}\|\tilde{\sigma}\|^2 = - \int_\Omega (\tilde\sigma\nabla\cdot\bu)\tilde\sigma\mathrm{d}\mathbf{x} - 2\int_\Omega (\bu\cdot\nabla\hat\sigma)\tilde\sigma\mathrm{d}\mathbf{x} - \int_\Omega (\tilde\sigma\nabla\cdot\bu)\hat\sigma\mathrm{d}\mathbf{x} -2\int_\Omega (\bu\cdot\nabla\tilde\sigma)\tilde\sigma\mathrm{d}\mathbf{x}.
\end{align}
Taking $L^2$ inner product of \eqref{aa5} with $\bu$, we get 
\begin{align}\label{AA31}
 \frac{\mathrm{d}}{\mathrm{d}t}\|\bu\|^2 + 2\alpha\|\bu\|^2 =  &\,- \int_\Omega (\bu\cdot\nabla\tilde\sigma)\tilde\sigma \mathrm{d}\mathbf{x} - \int_\Omega  (\bu\cdot\nabla\hat\sigma)\tilde\sigma\mathrm{d}\mathbf{x} - \int_\Omega (\bu\cdot\nabla\tilde\sigma)\hat\sigma\mathrm{d}\mathbf{x} \notag\\
&\, -2\int_\Omega [(\bu\cdot\nabla)\bu]\cdot\bu \mathrm{d}\mathbf{x} + 2\beta\int_\Omega \bu\cdot \nabla\tilde\phi \mathrm{d}\mathbf{x}.
\end{align}
By adding \eqref{AA30} and \eqref{AA31}, we obtain 
\begin{align}\label{AA32}
\frac{\mathrm{d}}{\mathrm{d}t}\|(\tilde{\sigma},\bu)\|^2 +2 \alpha\|\bu\|^2 = &\, -\int_\Omega [\nabla\cdot(\tilde\sigma\bu)]\tilde\sigma \mathrm{d}\mathbf{x} - 3\int_\Omega  (\bu\cdot\nabla\hat\sigma)\tilde\sigma\mathrm{d}\mathbf{x} - \int_\Omega [\nabla\cdot(\tilde\sigma\bu)]\hat\sigma \mathrm{d}\mathbf{x}\notag\\
&\, -\int_\Omega \bu\cdot\nabla (\tilde\sigma^2)\mathrm{d}\mathbf{x} - \int_\Omega (\bu\cdot\nabla) |\bu|^2\mathrm{d}\mathbf{x} + 2\beta\int_\Omega \bu\cdot \nabla\tilde\phi \mathrm{d}\mathbf{x}.
\end{align}
After integrating by parts, we update \eqref{AA32}
\begin{align}\label{AA33}
\frac{\mathrm{d}}{\mathrm{d}t}\|(\tilde{\sigma},\bu)\|^2 +2 \alpha\|\bu\|^2 = &\, \int_\Omega \tilde\sigma\bu\cdot\nabla\tilde\sigma \mathrm{d}\mathbf{x} - 2\int_\Omega  (\bu\cdot\nabla\hat\sigma)\tilde\sigma\mathrm{d}\mathbf{x} + \int_\Omega (\nabla\cdot\bu)\tilde\sigma^2\mathrm{d}\mathbf{x} \notag\\
&\, + \int_\Omega (\nabla\cdot\bu) |\bu|^2\mathrm{d}\mathbf{x} + 2\beta\int_\Omega \bu\cdot \nabla\tilde\phi \mathrm{d}\mathbf{x} \ \triangleq \ R_1+\cdots + R_5.
\end{align}
We emphasize that the derivation of \eqref{AA33} is a systematic approach that can be applied to any order temporal derivative of the solution.

For the nonlinear terms on the right-hand side of\eqref{AA33}, by H\"older's inequality, we have 
\begin{align}\label{AA34}
|R_1|+|R_2|+|R_3|+|R_4| &\leqslant \big(\|\nabla\tilde\sigma\|_{L^\infty} + 2\|\nabla\hat\sigma\|_{L^\infty}\big) \|\tilde\sigma\| \|\bu\| + \|\nabla\cdot\bu\|_{L^\infty} \big(\|\tilde\sigma\|^2+\|\bu\|^2\big) \notag \\
&\leqslant \big(\|\nabla\tilde\sigma\|_{L^\infty} +\|\nabla\cdot\bu\|_{L^\infty} + \|\nabla\hat\sigma\|_{L^\infty}\big) \big(\|\tilde\sigma\|^2+\|\bu\|^2\big).
\end{align}
Using the Sobolev embedding $H^2 \hookrightarrow L^\infty$ and the definition of $\mathtt{E}(t)$, we upgrade \eqref{AA34} as 
\begin{align}\label{AA35}
|R_1|+|R_2|+|R_3|+|R_4| \leqslant \big(\mathtt{E}^\frac12+ \|\nabla\hat\sigma\|_{L^\infty}\big) \mathtt{E}.
\end{align}
For $R_5$, by Cauchy's inequality, we have
\begin{align}\label{AA36}
|R_5| \leqslant \alpha\|\bu\|^2 + \alpha^{-1}\beta^2\|\nabla\tilde\phi\|^2 \leqslant \alpha\|\bu\|^2 + Cd^{-1}\|\tilde\sigma\|^2,
\end{align}
where \eqref{AA9a} is applied. Substituting \eqref{AA35}--\eqref{AA36} into \eqref{AA33} yields 
\begin{align}\label{AA37}
\frac{\mathrm{d}}{\mathrm{d}t}\|(\tilde{\sigma},\bu)\|^2 + \alpha\|\bu\|^2 \leqslant \big(\mathtt{E}^\frac12+ \|\nabla\hat\sigma\|_{L^\infty} + d^{-1}\big) \mathtt{E}.
\end{align}
In a similar fashion, it can be shown that 
\begin{align}\label{AA38}
\frac{\mathrm{d}}{\mathrm{d}t}\Big(\sum_{k=0}^3\|\partial_t^k(\tilde{\sigma},\bu)\|^2\Big) + \alpha \sum_{k=0}^3\|\partial_t^k\bu\|^2 \leqslant C\big(\mathtt{E}^\frac12+ \|\nabla\hat\sigma\|_{L^\infty} + d^{-1}\big) \mathtt{E}.
\end{align}
This completes the third step of the full energy scheme. Next, we recover the dissipation of $\tilde\sigma$.

\vspace{.1 in}

{\bf Recovery of density dissipation}: We first derive a wave-type equation for $\tilde\sigma$. Taking temporal derivative of \eqref{aa4} and spatial divergence of \eqref{aa5}, we get
\begin{subequations}\label{AA39}
\begin{alignat}{2}
2 \tilde{\sigma}_{tt} + (\hat{\sigma} \nabla\cdot\bu)_t &=-\big[\bu\cdot \nabla{\tilde{\sigma}}_t + \nabla\cdot(\tilde{\sigma} \bu)_t +\bu_t\cdot\nabla\tilde\sigma + 2 \bu_t\cdot \nabla\hat{\sigma}\big], \label{aa9}\\
2 (\nabla\cdot\bu)_t + \nabla\cdot(\hat{\sigma}\nabla{\tilde{\sigma}}) &=-\nabla\cdot \big[2 (\bu\cdot \nabla )\bu + \tilde{\sigma}\nabla{\tilde{\sigma}} + \tilde{\sigma} \nabla{\hat{\sigma}} + 2 \alpha \bu - 2 \beta \nabla{\tilde{\phi}}\big]. \label{aa10}
\end{alignat}
\end{subequations}
To simplify the presentation, we define 
\begin{align}
\mathsf{M} &\triangleq \bu\cdot \nabla{\tilde{\sigma}}_t + \nabla\cdot(\tilde{\sigma} \bu)_t +\bu_t\cdot\nabla\tilde\sigma + 2 \bu_t\cdot \nabla\hat{\sigma}, \label{AA40}\\
{\boldsymbol{\mathsf{N}}} &\triangleq 2 (\bu\cdot \nabla )\bu + \tilde{\sigma}\nabla{\tilde{\sigma}} + \tilde{\sigma} \nabla{\hat{\sigma}} + 2 \alpha \bu - 2 \beta \nabla{\tilde{\phi}}. \label{AA41}
\end{align}
Multiplying \eqref{aa10} by $\hat\sigma$, then subtracting the result from \eqref{aa9}$\times 2$, we have 
\begin{align}\label{AA42}
4\tilde{\sigma}_{tt} - \nabla\cdot (\hat\sigma^2\nabla\tilde\sigma) = -2\mathsf{M} + \nabla\cdot(\hat\sigma {\boldsymbol{\mathsf{N}}}) - {\boldsymbol{\mathsf{N}}}\cdot \nabla\hat\sigma -\hat\sigma\nabla\tilde\sigma\cdot\nabla\hat\sigma,
\end{align}
which constitutes a non-homogeneous wave equation for $\tilde\sigma$ with non-constant speed. Taking $L^2$ inner product of \eqref{AA42} with $-\tilde\sigma$, we can show that
\begin{align}\label{AA43}
-4\frac{\mathrm{d}}{\mathrm{d}t} \langle\tilde\sigma_t,\tilde\sigma\rangle + 4\|\tilde\sigma_t\|^2 = &\,2\int_\Omega \tilde\sigma \mathsf{M} \mathrm{d}\mathbf{x} - \int_\Omega \tilde\sigma\nabla\cdot[\hat\sigma({\boldsymbol{\mathsf{N}}}+\hat\sigma\nabla\tilde\sigma)] \mathrm{d}\mathbf{x} + \int_\Omega \tilde\sigma {\boldsymbol{\mathsf{N}}} \cdot \nabla\hat\sigma \mathrm{d}\mathbf{x}\notag\\
&\,+ \int_\Omega \tilde\sigma \hat\sigma\nabla\tilde\sigma\cdot\nabla\hat\sigma \mathrm{d}\mathbf{x},
\end{align}
where $\langle\cdot,\cdot\rangle$ denotes the standard inner product in $L^2$. For the first integral on the right-hand side of\eqref{AA43}, after integrating by parts, we have
\begin{align}\label{AA44}
2\int_\Omega \tilde\sigma \mathsf{M} \mathrm{d}\mathbf{x} = -2\int_\Omega\big[ \tilde\sigma_t\nabla\cdot(\tilde\sigma\bu) +(\tilde\sigma\bu)_t\cdot\nabla\tilde\sigma - \tilde\sigma\bu_t\cdot\nabla\tilde\sigma - 2\tilde\sigma\bu_t\cdot\nabla\hat\sigma\big]\mathrm{d}\mathbf{x}.
\end{align}
It follows from H\"older's inequality and Sobolev embedding that 
\begin{align}\label{AA45}
\Big|2\int_\Omega \tilde\sigma \mathsf{M} \mathrm{d}\mathbf{x}\Big| \leqslant C\big(\mathtt{E}^\frac12+ \|\nabla\hat\sigma\|_{L^\infty}\big) \mathtt{E}.
\end{align}
For the second integral on the right-hand side of\eqref{AA43}, according to \eqref{aa5} and the boundary condition for $\bu$, it holds that $({\boldsymbol{\mathsf{N}}}+\hat\sigma\nabla\tilde\sigma) \cdot\mathbf{n}|_{\partial\Omega} =0$. Hence, after integrating by parts, we obtain
\begin{align}\label{AA46}
- \int_\Omega \tilde\sigma\nabla\cdot[\hat\sigma({\boldsymbol{\mathsf{N}}}+\hat\sigma\nabla\tilde\sigma)] \mathrm{d}\mathbf{x} = \int_\Omega \hat\sigma({\boldsymbol{\mathsf{N}}}+\hat\sigma\nabla\tilde\sigma)\cdot\nabla \tilde\sigma \mathrm{d}\mathbf{x} = -2\int_\Omega \hat\sigma \bu_t \cdot\nabla \tilde\sigma \mathrm{d}\mathbf{x}.
\end{align}
This implies that 
\begin{align}\label{AA47}
\Big|-\int_\Omega \tilde\sigma\nabla\cdot[\hat\sigma({\boldsymbol{\mathsf{N}}}+\hat\sigma\nabla\tilde\sigma)] \mathrm{d}\mathbf{x}\Big| \leqslant \|\hat\sigma\|_{L^\infty}\big(\|\bu_t\|^2 +\|\nabla\tilde\sigma\|^2 \big).
\end{align}
Further adopting \eqref{AA18}, we update \eqref{AA47} as  
\begin{align}\label{AA48}
\Big|-\int_\Omega \tilde\sigma\nabla\cdot[\hat\sigma({\boldsymbol{\mathsf{N}}}+\hat\sigma\nabla\tilde\sigma)] \mathrm{d}\mathbf{x}\Big| \leqslant C\big(\|(\bu_t,\bu)\|^2 + \mathtt{E}^2\big),
\end{align}
where $\mathtt{E}$ is defined by \eqref{AA5}. For the last two integrals on the right-hand side of\eqref{AA43}, it is not difficult to see that they are bounded similarly as \eqref{AA45}. Collectively, we obtain 
\begin{align}\label{AA49}
-4\frac{\mathrm{d}}{\mathrm{d}t} \langle\tilde\sigma_t,\tilde\sigma\rangle + 4\|\tilde\sigma_t\|^2 \leqslant C\big(\mathtt{E}^\frac12+ \|\nabla\hat\sigma\|_{L^\infty}\big) \mathtt{E} + C\|(\bu_t,\bu)\|^2.
\end{align}
In a similar fashion, by taking temporal derivatives of \eqref{AA42}, it can be shown that 
\begin{align}\label{AA50}
-4\frac{\mathrm{d}}{\mathrm{d}t} \Big(\sum_{k=1}^3 \langle\partial_t^k \tilde\sigma,\partial_t^{k-1}\tilde\sigma\rangle\Big) + 4\sum_{k=1}^3\|\partial_t^k\tilde\sigma\|^2 \leqslant C\big(\mathtt{E}^\frac12+ \|\nabla\hat\sigma\|_{L^\infty}\big) \mathtt{E} + C\sum_{k=0}^3\|\partial_t^k\bu\|^2.
\end{align}
This completes the recovery of density dissipation.

\vspace{.1 in}

{\bf Estimation of curl}: The energy estimates derived in $\S$3.3 and $\S$3.4 are related to the energetic quantity $\mathtt{E}$ defined in \eqref{AA5}. For the other, i.e., $\mathtt{V}_1 =\|\omega\|_2^2$, by taking the curl of \eqref{aa2}, we have
\begin{align}\label{AA51}
\omega_t + (\bu\cdot\nabla) \omega + \omega\nabla\cdot\bu = -\alpha\omega.
\end{align}
By taking mixed spatial-temporal derivatives of up to the second order and applying integration-by-parts, it can be shown that 
\begin{align}\label{AA52}
\frac{\mathrm{d}}{\mathrm{d}t} \mathtt{V}_1 + \alpha \mathtt{V}_1 \leqslant C\mathtt{E}^{\frac32}.
\end{align}
At this stage, we have completed the estimates for each component defined in \eqref{AA5}. The next step involves coupling these components to close the full energy scheme.

\vspace{.1 in}

{\bf Closing of energy estimates}: For reader's convenience, we list below the energetic inequalities \eqref{AA38} and \eqref{AA50}: 
\begin{align}
\frac{\mathrm{d}}{\mathrm{d}t}\Big(\sum_{k=0}^3\|\partial_t^k(\tilde{\sigma},\bu)\|^2\Big) + \alpha \sum_{k=0}^3\|\partial_t^k\bu\|^2 &\leqslant C\big(\mathtt{E}^\frac12+ \|\nabla\hat\sigma\|_{L^\infty} + d^{-1}\big) \mathtt{E},\label{AA53}\\
-4\frac{\mathrm{d}}{\mathrm{d}t} \Big(\sum_{k=1}^3 \langle\partial_t^k \tilde\sigma,\partial_t^{k-1}\tilde\sigma\rangle\Big) + 4\sum_{k=1}^3\|\partial_t^k\tilde\sigma\|^2 &\leqslant C\big(\mathtt{E}^\frac12+ \|\nabla\hat\sigma\|_{L^\infty}\big) \mathtt{E} + \mathsf{C}_1\sum_{k=0}^3\|\partial_t^k\bu\|^2,\label{AA54} 
\end{align}
where, for clarity, we have relabeled the two constants appearing on the right-hand sides of the inequalities. Multiplying \eqref{AA53} by a suitably large constant $\mathsf{L}_1$ such that $\mathsf{L}_1>\max\{4,\alpha ^{-1} \mathsf{C}_1\}$, then adding the result to \eqref{AA54}, we have 
\begin{align}\label{AA55}
\frac{\mathrm{d}}{\mathrm{d}t} \mathtt{X}_1 + \mathtt{Y}_1 \leqslant C\big(\mathtt{E}^\frac12+ \|\nabla\hat\sigma\|_{L^\infty} + d^{-1}\big) \mathtt{E},
\end{align}
where the quantities on the left are defined by 
\begin{align}
\mathtt{X}_1 &\triangleq \mathsf{L}_1 \sum_{k=0}^3\|\partial_t^k(\tilde{\sigma},\bu)\|^2 - 4 \sum_{k=1}^3\int_\Omega \big(\partial_t^k \tilde\sigma\big) \big(\partial_t^{k-1}\tilde\sigma\big) \mathrm{d}\mathbf{x}, \label{AA56}\\
\mathtt{Y}_1 &\triangleq \mathsf{L}_1\alpha \sum_{k=0}^3\|\partial_t^k\bu\|^2 - \mathsf{C}_1\sum_{k=0}^3\|\partial_t^k\bu\|^2 + 4\sum_{k=1}^3\|\partial_t^k\tilde\sigma\|^2. \label{AA57}
\end{align}
By virtue of the choice of $\mathsf{L}_1$, we know that 
\begin{align}
\mathtt{X}_1 &\cong \sum_{k=0}^3\|\partial_t^k(\tilde{\sigma},\bu)\|^2= \mathtt{U}_1+\|\tilde\sigma\|^2, \label{AA58}\\
\mathtt{Y}_2 &\cong \sum_{k=0}^3\|\partial_t^k\bu\|^2 + \sum_{k=1}^3\|\partial_t^k\tilde\sigma\|^2 =\mathtt{U}. \label{AA59}
\end{align}
Taking the sum of \eqref{AA55} and \eqref{AA52} yields
\begin{align}\label{AA63}
\frac{\mathrm{d}}{\mathrm{d}t} (\mathtt{X}_1 +\mathtt{V}_1) + \mathtt{Y}_1 + \alpha\mathtt{V}_1 \leqslant C\big(\mathtt{E}^\frac12+ \|\nabla\hat\sigma\|_{L^\infty} + d^{-1}\big) \mathtt{E}.
\end{align}
Using the qualitative relationships \eqref{AA58}, \eqref{AA59}, \eqref{AA18}, and \eqref{AA29}, we upgrade \eqref{AA63} as 
\begin{align}\label{AA64}
\frac{\mathrm{d}}{\mathrm{d}t} (\mathtt{X}_1 +\mathtt{V}_1) + C(\mathtt{X}_1 + \mathtt{V}_1) \leqslant C\big(\mathtt{E}^\frac12+ \|\nabla\hat\sigma\|_{L^\infty} + d^{-1}\big) (\mathtt{X}_1 + \mathtt{V}_1).
\end{align}
Since $\|\nabla\hat\sigma\|_{L^\infty} \to 0$ as $d\to\infty$, when $\tilde\delta$ is sufficiently small and $d$ is sufficiently large, we obtain
\begin{align}\label{AA65}
\frac{\mathrm{d}}{\mathrm{d}t} (\mathtt{X}_1 +\mathtt{V}_1) + C(\mathtt{X}_1 + \mathtt{V}_1) \leqslant 0,
\end{align}
which implies the uniform $L^\infty_t$ and $L^1_t$ estimates and exponential decay of $\mathtt{X}_1(t) +\mathtt{V}_1(t)$. Since $\mathtt{E}(t) \lesssim \mathtt{X}_1(t)+\mathtt{V}_1(t)$, we conclude that $\mathtt{E}(t)$ is bounded in $(L^\infty\cap L^1) (0,\infty)$ and decays exponentially to zero as $t\to\infty$, which, along with elliptic regularity estimates, implies the same property for $\sum_{k=0}^3\|(\partial_t^k\tilde\phi)(t)\|_{H^{5-k}}^2$. This completes the proof of Theorem \ref{main}.


\section{Steady state and its stability under \eqref{bc1} $\&$ \eqref{bc3}}

In this section, we investigate the steady-state solution for \eqref{eqn:vasculogenesismodel} with the ``$\gamma$-law pressure":
\begin{subequations}\label{HP}
\begin{alignat}{3}
\rho_t + \nabla \cdot(\rho \mathbf{u}) &=0, \label{HP1}\\
(\rho \mathbf{u})_t + \nabla \cdot(\rho \mathbf{u}\otimes \mathbf{u})+\nabla (A_0\rho^\gamma) &=-\alpha\rho \mathbf{u}+\beta\rho\nabla\phi, \label{HP2}\\
\tau\phi_t - d\Delta\phi + a\phi &=b\rho, \label{HP3}
\end{alignat}
\end{subequations}
under the boundary conditions \eqref{bc1} and \eqref{bc3}, where $\gamma \geqslant 1$ is a constant.

\subsection{Finding steady state} 

We emphasize that our primary objective is to identify asymptotically stable steady-state solutions. As shown in the previous section, the Poincaré inequality for the perturbed chemical concentration is indispensable in our analysis--indeed, without it, the entire energy framework would fail to close. In the Dirichlet boundary condition case, both the time-evolutionary and steady-state concentrations satisfy the inequality for all time. However, under the Neumann boundary condition for $\phi$, the Poincaré inequality only holds if the perturbed concentration is spatially mean-free. Yet, the spatial mean of $\phi$ is not a steady profile; rather, it evolves dynamically. To verify this, we integrate \eqref{HP3} over $\Omega$ to obtain:
\begin{align}\label{BB1}
\tau\frac{\mathrm{d}}{\mathrm{d}t} \int_\Omega \phi(\mathbf{x},t)\mathrm{d}\mathbf{x} = -a \int_\Omega \phi(\mathbf{x},t)\mathrm{d}\mathbf{x} + b\int_\Omega \rho(\mathbf{x},t)\mathrm{d}\mathbf{x}.
\end{align}
Let $\hat\phi$ and $\hat\rho$ denote the spatial average of $\phi$ and $\rho$, respectively. Under the no-normal-flow boundary condition, $\hat\rho=\hat{\rho}_0$ remains constant. From the equation above, we derive:
\begin{align}\label{BB2}
\hat\phi(t) = \frac{b}{a}\hat\rho_0 + \Big(\hat\phi_0-\frac{b}{a}\hat\rho_0\Big) e^{-\frac{a}{\tau} t},
\end{align}
where $\hat\phi_0$ is the initial spatial average of $\phi$. This result necessitates abandoning the conventional approach of subtracting a steady-state solution from the time-dependent solution to establish asymptotic decay of perturbations. Nevertheless, \eqref{BB2} reveals that if $\phi-\hat\phi$ vanishes, the steady-state value of $\phi$ must be $\frac{b}{a}\hat\rho_0$. Consequently, the asymptotically stable steady-state solution of \eqref{HP} under boundary conditions \eqref{bc1} and \eqref{bc3} should take the form:
\begin{align}\label{BB3}
(\rho_\infty,\bu_\infty,\phi_\infty) = \Big(\hat\rho_0,\mathbf{0},\frac{b}{a}\hat\rho_0\Big).
\end{align}
Next, we demonstrate that this steady state is locally asymptotically stable, provided the initial data and diffusion coefficient $d$ satisfy conditions analogous to those in Theorem \ref{main}.

\subsection{Stability of steady state}

Our stability results are summarized in the following theorem.

\begin{theorem}\label{main1}
Consider the initial-boundary value problem of \eqref{HP} with initial data $(\rho_0,\bu_0,\phi_0)$ and boundary conditions \eqref{bc1} $\&$ \eqref{bc3}. Let $A_0,\alpha,\beta,\tau,d,a,b>0$ and $\gamma\geqslant 1$ be constants, and let $\hat\rho_0$ and $\hat\phi_0$ be the spatial average of $\rho_0$ and $\phi_0$, respectively. Suppose that there is a large constant $d_1>0$ such that $d \geqslant d_1$, and that the initial data satisfy the following conditions:
\begin{itemize}
\item $(\rho_0,\bu_0,\phi_0)$ are compatible with the boundary conditions;
\item $\rho_0>0$, $(\rho_0,\bu_0) \in [H^3(\Omega)]^3$, and $\phi_0\in H^5(\Omega)$;
\item there is a small constant $\delta_1>0$ such that $\|(\rho_0-\hat{\rho}_0,\bu_0)\|^2_{H^3}+\|\phi_0-\hat\phi_0\|_{H^5}^2\leqslant \delta_1$.
\end{itemize}
Then there exists a unique solution to the initial-boundary value problem, such that for $\forall t>0$,
\begin{equation*}
\begin{aligned}
\mathfrak{E}(t) + \int_0^t \big(\mathfrak{E}(\tau) +\|\tilde\phi_{tt}(\tau)\|_{H^2}^2\big)\mathrm{d}\tau \leqslant C,
\end{aligned}
\end{equation*}
where $\mathfrak{E}(t) \triangleq \|(\tilde{\rho}, \bu)(t)\|^2_3 + \sum_{k=0}^2 \|(\partial _t^k \tilde{\phi})(t)\|^2_{H^{5-2k}}$, $\tilde{\rho}=\rho-\hat{\rho}_0$, $\tilde{\phi}=\phi-\hat\phi(t)$, and $\hat\phi(t)$ is given by \eqref{BB2}. Moreover, there are positive constants $\eta_3$ and $\eta_4$, such that for $\forall t>0$, 
$\mathfrak{E}(t) \leqslant \eta_3 e^{-\eta_4 t}$.
\end{theorem}

\begin{rem}
In view of \eqref{BB2}, we see that $\hat\phi(t) \to \frac{b}{a}\hat\rho_0$ exponentially in time. As a consequence of Theorem \ref{main1}, the steady state specified in \eqref{BB3} is locally exponentially stable. 
\end{rem}

\begin{rem}
The proof of Theorem \ref{main1} adopts a strategy analogous to that employed in the preceding section. The key distinction arises in the symmetrization of the fluid component: for $\gamma>1$, we utilize the general sound-speed transformation \eqref{SST}, whereas the case $\gamma=1$ requires application of a diagonal symmetrizer. Due to this dichotomy, the proof is naturally divided into two subsections.
\end{rem}

\begin{rem}
The analysis in this section is not restricted to the square domain case. For obtuse-angled polygons or general bounded domains with smooth boundaries, one can utilize coordinate transformation and local boundary straightening (c.f.\,\cite{MN}), and adapt the proof techniques in this section to establish similar results as in Theorem \ref{main1}. 
\end{rem}

\subsection{Proof of Theorem \ref{main1} for $\gamma > 1$} Applying the general sound-speed transformation:
\begin{align}\label{SST}
\sigma=\mathfrak{A}\rho^\kappa,\quad \text{where}\quad \mathfrak{A}=\frac{2\sqrt{\gamma A_0}}{\gamma-1},\quad \kappa=\frac{\gamma-1}{2},
\end{align}
we obtain the following system in terms of the perturbed variables:
\begin{subequations}\label{N1}
\begin{alignat}{3}
\tilde{\sigma}_t + \bu\cdot \nabla{\tilde{\sigma}} + \kappa \tilde{\sigma} \nabla\cdot \bu + \kappa \hat{\sigma} \nabla\cdot\bu &=0, \label{N1a}\\
\bu_t + (\bu\cdot \nabla )\bu + \kappa \tilde{\sigma}\nabla{\tilde{\sigma}} + \kappa \hat{\sigma}\nabla{\tilde{\sigma}}  
&= - \alpha \bu +  \beta \nabla{\tilde{\phi}}, \label{N1b}\\
\tau \tilde\phi_t - d \Delta\tilde{\phi} + a \tilde{\phi} &=  b\mathfrak{A}^{-\frac{1}{\kappa}} \big[(\tilde{\sigma} +  \hat{\sigma})^{\frac{1}{\kappa}} - \hat{\sigma}^{\frac{1}{\kappa}}\big], \label{N1c}
\end{alignat}
\end{subequations}
where $\tilde\sigma=\sigma-\hat\sigma$, $\hat\sigma=\mathfrak{A}\hat\rho_0^\kappa$, and $\hat\rho_0$ and $\tilde\phi$ are defined in Theorem \ref{main1}. The initial and boundary conditions for system \eqref{N1} read as
\begin{align*}
\big(\tilde{\sigma}, \bu,\tilde\phi)(\mathbf{x},0\big) &= \big(\mathfrak{A}\rho_0^\kappa - \mathfrak{A}\hat\rho_0^\kappa, \bu_0,\phi_0-\hat\phi_0\big)(\mathbf{x}), \\
\bu\cdot\mathbf{n}|_{\partial \Omega}&=0=\nabla\tilde{\phi}\cdot\mathbf{n}|_{\partial \Omega}. 
\end{align*}
In the subsequent analysis, we shall carry out energy estimates under the {\it a priori} assumption:
\begin{align}\label{BB4}
\esssup_{t\in [0,T]} \mathfrak{E}(t) \leqslant \delta,
\end{align}
where $\mathfrak{E}$ is defined in Theorem \ref{main1}. Here, without abuse of notation, we reuse the symbols $T$ and $\delta$ from \eqref{AA4} to denote the lifespan of the local solution and a small positive constant, respectively. Next, we apply the same strategy as in the previous section to prove Theorem \ref{main1}.

\vspace{.1 in}

{\bf Reduction of total energy}: To reduce the total energy $\mathfrak{E}(t)$, we define the following:
\begin{align}\label{BB5}
\mathtt{U}_2(t) \triangleq \|(\bu_{tt},\tilde\sigma_{tt})(t)\|^2_{H^1},\quad \mathtt{V}_2(t) \triangleq \|\omega(t)\|_{H^2}^2+\|\omega_t\|_{H^1}^2,\quad \mathtt{W}(t) \triangleq \|(\nabla\tilde\phi,\nabla\tilde\phi_t,\nabla\tilde\phi_{tt})(t)\|^2.
\end{align}
We shall show that $\mathfrak{E}(t)$ is controlled by the quantities specified in \eqref{BB5}.

{\bf Step 1.} First of all, note that since $\tilde\rho$ is mean-free, $\tilde\sigma$ satisfies Poincar\'e's inequality when $\delta$ is sufficiently small. In addition, since $\tilde\phi$ is mean free, it also satisfies the Poincar\'e inequality.

{\bf Step 2.} Using \eqref{N1a}--\eqref{N1b}, we deduce that under \eqref{BB4},
 \begin{align}
\|\bu\|_{H^3}^2 &\leqslant C\big(\|\tilde\sigma_t\|_{H^2}^2 + \|\omega\|_{H^2}^2+ \mathfrak{E}^2\big), \label{BB6} \\
\|\tilde\sigma\|_{H^3}^2 &\leqslant C\big(\|(\bu_t,\bu)\|_{H^2}^2+\|\nabla\tilde\phi\|_{H^2}^2+\mathfrak{E}^2\big). \label{BB7}
\end{align}
Here, we note that $\tilde\phi$ appears on the right-hand side of\eqref{BB7} because $\tilde\phi$ satisfies a parabolic equation rather than an elliptic one. Similarly, by reduction, we can show that 
 \begin{align}
\|\bu_t\|_{H^2}^2 &\leqslant C\big(\|\tilde\sigma_{tt}\|_{H^1}^2 + \|\omega_t\|_{H^1}^2+ \mathfrak{E}^2\big), \label{BB8} \\
\|\tilde\sigma_t\|_{H^2}^2 &\leqslant C\big(\|(\bu_{tt},\bu_t)\|_{H^1}^2+\|\nabla\tilde\phi_t\|_{H^1}^2+\mathfrak{E}^2\big). \label{BB9}
\end{align}
Combining \eqref{BB6}--\eqref{BB9} yields
 \begin{align}\label{BB10}
\|(\bu,\tilde\sigma)\|_{H^3}^2 \leqslant C\big(\|(\bu_{tt},\tilde\sigma_{tt})\|_{H^1}^2 + \|\omega\|_{H^2}^2 + \|\omega_t\|_{H^1}^2 + \|\nabla\tilde\phi\|_{H^2}^2 + \|\nabla\tilde\phi_t\|_{H^1}^2 +\mathfrak{E}^2\big).
\end{align}
Since $\tilde\phi$ satisfies the zero Neumann boundary condition, by \eqref{AA14}, we have
\begin{align}\label{BB11}
\|\nabla\tilde\phi\|_{H^2}^2 +\|\nabla\tilde\phi_t\|_{H^1}^2 \leqslant C\big(\|\Delta\tilde\phi\|_{H^1}^2 +\|\Delta\tilde\phi_t\|^2 + \|\nabla\tilde\phi\|^2 + \|\nabla\tilde\phi_t\|^2\big).
\end{align}
Using \eqref{N1c}, the spatial derivatives can be reduced as follows:  
\begin{align}\label{BB12}
\|\Delta\tilde\phi\|_{H^1}^2 + \|\Delta\tilde\phi_t\|^2 \leqslant Cd^{-2}\big(\|\tilde\phi_t\|_{H^1}^2 + \|\tilde\phi_{tt}\|^2 + \|\tilde\phi\|_{H^1}^2 + \|\tilde\phi_{t}\|^2 + \mathfrak{E}\big).
\end{align}
Since $\tilde\phi$ and its temporal derivatives satisfy Poincar\'e's inequality, we obtain from \eqref{BB11}--\eqref{BB12}:
\begin{align}\label{BB13}
\|\nabla\tilde\phi\|_{H^2}^2 + \|\nabla\tilde\phi_t\|_{H^1}^2 \leqslant C\mathtt{W} + Cd^{-2}\mathfrak{E}.
\end{align}
Substituting \eqref{BB13} into \eqref{BB8}--\eqref{BB10} yields the quantitative relationship:
 \begin{align}\label{BB14}
\sum_{k=0}^2\|(\partial_t^k\bu,\partial_t^k\tilde\sigma)\|_{H^{3-k}}^2 \leqslant C\big(\|(\bu_{tt},\tilde\sigma_{tt})\|_{H^1}^2 + \|\omega\|_{H^2}^2 + \|\omega_t\|_{H^1}^2 + \mathtt{W} +d^{-2}\mathfrak{E}+\mathfrak{E}^2\big).
\end{align}
Furthermore, by \eqref{N1a}--\eqref{N1b}, we can show that 
\begin{align}\label{BB15}
\|(\bu_{ttt},\tilde\sigma_{ttt})\|^2 \leqslant C\big(\|(\bu_{tt},\tilde\sigma_{tt})\|_{H^1}^2 + \|\nabla\tilde\phi_{tt}\|^2+\mathfrak{E}^2\big).
\end{align}
Combining \eqref{BB14} and \eqref{BB15}, we arrive at 
\begin{align}\label{BB16}
\|(\tilde\sigma,\bu)(t)\|_3^2 \leqslant C\big(\mathtt{U}_2+\mathtt{V}_2+\mathtt{W}+d^{-2}\mathfrak{E}+\mathfrak{E}^2\big),
\end{align}
where we adopted the notations from \eqref{BB5}.

{\bf Step 3.} According to \eqref{N1c}, \eqref{AA14}, and Poincar\'e's inequality, we have
\begin{align}\label{BB17}
\|\tilde\phi\|_{H^5}^2 \leqslant C\big(\|\Delta\tilde\phi\|_{H^3}^2+\|\nabla\tilde\phi\|^2\big) \leqslant C\big(\|\tilde\phi_t\|_{H^3}^2+\|\tilde\phi\|_{H^3}^2\big)+Cd^{-2}\mathfrak{E}.
\end{align}
Similarly, it holds that
\begin{align}\label{BB18}
\|\tilde\phi_t\|_{H^3}^2 \leqslant C\big(\|\Delta\tilde\phi_t\|_{H^1}^2 + \|\nabla\tilde\phi_t\|^2\big) \leqslant C\big(\|\nabla\tilde\phi_{tt}\|^2+ \|\nabla\tilde\phi_t\|^2\big) + Cd^{-2}\mathfrak{E}.
\end{align}
Moreover, it can be shown that 
\begin{align}\label{BB19}
\|\tilde\phi\|_{H^3}^2 \leqslant C\big(\|\Delta\tilde\phi\|_{H^1}^2+\|\nabla\tilde\phi\|^2\big) \leqslant C\big(\|\nabla\tilde\phi_t\|^2+\|\nabla\tilde\phi\|^2\big)+Cd^{-2}\mathfrak{E}.
\end{align}
Combining \eqref{BB17}--\eqref{BB19}, we obtain 
\begin{align}\label{BB20}
\sum_{k=0}^2\|(\partial_t^k\tilde\phi)(t)\|_{H^{5-2k}}^2 \leqslant C\mathtt{W} + Cd^{-2}\mathfrak{E}.
\end{align}
When $\mathfrak{E}$ is small and $d$ is large, \eqref{BB16} and \eqref{BB20} imply
\begin{align}\label{BB21}
\mathfrak{E}(t) \leqslant C\big[\mathtt{U}_2(t)+\mathtt{V}_2(t)+\mathtt{W}(t)\big].
\end{align}
Next, we deal with the estimation of the quantities on the right-hand side of\eqref{BB21}.

\vspace{.1 in}

{\bf Estimation of reduced energy}: First of all, similar to \eqref{AA37}, we can show that 
\begin{align}\label{BB21a}
\frac{\mathrm{d}}{\mathrm{d}t}\|(\tilde{\sigma},\bu)\|^2 + \alpha\|\bu\|^2 \leqslant \mathfrak{E}^\frac32 + C\|\nabla\tilde\phi\|^2.
\end{align}
Next, taking $\partial_x$ of \eqref{N1a} and \eqref{N1b}, we have
\begin{align}
\tilde{\sigma}_{xt} +  \bu\cdot \nabla{\tilde{\sigma}}_x + \kappa \tilde{\sigma} \nabla\cdot \bu_x + \kappa\hat{\sigma} \nabla\cdot\bu_x &=- \bu_x\cdot \nabla{\tilde{\sigma}} - \kappa \tilde{\sigma}_x \nabla\cdot \bu, \label{BB22}\\
\bu_{xt} + (\bu\cdot \nabla )\bu_x + \kappa\tilde{\sigma}\nabla{\tilde{\sigma}}_x  + \kappa\hat{\sigma}\nabla{\tilde{\sigma}}_x &=- \alpha \bu_x +  \beta \nabla\tilde{\phi}_x -  (\bu_x\cdot \nabla )\bu - \kappa\tilde{\sigma}_x\nabla{\tilde{\sigma}}, \label{BB23}
\end{align}
where we kept the terms requiring integration-by-parts on the left-hand sides of the equations. Taking $L^2$ inner product of \eqref{BB22} with $\tilde\sigma_x$, \eqref{BB23} with $\bu_x$, and performing standard energy estimates, we can show that
\begin{align}\label{BB24}
\frac{\mathrm{d}}{\mathrm{d}t} \|(\tilde\sigma_x,\bu_x)\|^2 + \alpha\|\bu_x\|^2 + 2\kappa\int_{\partial\Omega} (\tilde\sigma+\hat\sigma) \tilde\sigma_x\bu_x\cdot\mathbf{n}\mathrm{d}\mathbf{S} \leqslant  C\mathfrak{E}^\frac32 + C\|\nabla\tilde\phi_x\|^2.
\end{align}
Note that on $\partial\Omega$, because of the boundary conditions, there hold that 
\begin{align}
&v_x=0,\quad \tilde\sigma_y=0,\quad y=0,1,\ 0\leqslant x \leqslant 1, \label{BB25}\\
&u_y=0,\quad \tilde\sigma_x=0,\quad x=0,1,\ 0\leqslant y \leqslant 1. \label{BB26}
\end{align}
These imply the following:
\begin{equation}\label{BB27}
\tilde\sigma_x\bu_x\cdot\mathbf{n} =\left\{\begin{aligned} 
&\pm \tilde\sigma_x v_x =0,\quad y=0,1,\ 0\leqslant x \leqslant 1,\\
&\pm \tilde\sigma_x u_x =0,\quad x=0,1,\ 0\leqslant y \leqslant 1.
\end{aligned}
\right.
\end{equation}
Using \eqref{BB27}, we update \eqref{BB24} as 
\begin{align}\label{BB28}
\frac{\mathrm{d}}{\mathrm{d}t} \|(\tilde\sigma_x,\bu_x)\|^2 + \alpha\|\bu_x\|^2  \leqslant  C\mathfrak{E}^\frac32 + C\|\nabla\tilde\phi_x\|^2.
\end{align}
In a similar fashion, we have the following for the $y$-derivative: 
\begin{align}\label{BB29}
\frac{\mathrm{d}}{\mathrm{d}t} \|(\tilde\sigma_y,\bu_y)\|^2 + \alpha\|\bu_y\|^2  \leqslant C\mathfrak{E}^\frac32 + C\|\nabla\tilde\phi_y\|^2.
\end{align}
Combining \eqref{BB21a}, \eqref{BB28}, and \eqref{BB29}, we obtain
\begin{align}\label{BB30}
\frac{\mathrm{d}}{\mathrm{d}t} \|(\tilde\sigma,\bu)\|_{H^1}^2 + \alpha\|\bu\|_{H^1}^2  \leqslant C\mathfrak{E}^\frac32 + C\|\nabla\tilde\phi\|_{H^1}^2.
\end{align}
In a similar fashion, by taking temporal derivatives, it can be shown that 
\begin{align}\label{BB31}
\frac{\mathrm{d}}{\mathrm{d}t} \Big(\sum_{k=0}^2\|(\partial_t^k\tilde\sigma,\partial_t^k\bu)\|_{H^1}^2\Big) + \alpha\sum_{k=0}^2\|\partial_t^k\bu\|_{H^1}^2  \leqslant C\mathfrak{E}^\frac32 + C\sum_{k=0}^2\|\partial_t^k\nabla\tilde\phi\|_{H^1}^2.
\end{align}
As before, the dissipation of $\tilde\sigma$ is missing on the left-hand side of \eqref{BB31}, which is recovered next.

\vspace{.1 in}

{\bf Recovery of density dissipation}: First, similar to \eqref{AA49}, we can show that 
\begin{align}\label{BB32}
-\frac{\mathrm{d}}{\mathrm{d}t} \langle \tilde\sigma_t,\tilde\sigma\rangle + \|\tilde\sigma_t\|^2 \leqslant C\mathfrak{E}^\frac32 + \frac{\kappa\hat\sigma}{2}\|(\bu_t,\nabla\tilde\sigma)\|^2.
\end{align}
Next, similar to \eqref{AA42}, we have 
\begin{align}\label{BB33}
\tilde{\sigma}_{tt} = -\overline{\mathsf{M}} + \kappa\hat\sigma \nabla\cdot(\kappa\hat\sigma\nabla\tilde\sigma +\overline{{\boldsymbol{\mathsf{N}}}}),
\end{align}
where the quantities $\overline{\mathsf{M}}$ and $\overline{{\boldsymbol{\mathsf{N}}}}$ are defined by 
\begin{align*}
\overline{\mathsf{M}} &\triangleq \bu\cdot \nabla{\tilde{\sigma}}_t + \kappa\tilde\sigma\nabla\cdot \bu_t + \kappa \tilde\sigma_t\nabla\cdot\bu + \bu_t\cdot\nabla\tilde\sigma,\\
\overline{{\boldsymbol{\mathsf{N}}}} &\triangleq (\bu\cdot \nabla )\bu + \kappa \tilde{\sigma}\nabla{\tilde{\sigma}} +  \alpha \bu - \beta \nabla{\tilde{\phi}}.
\end{align*}
Taking $\partial_x$ of \eqref{BB33}, we obtain
\begin{align}\label{BB34}
\tilde{\sigma}_{xtt} = -\overline{\mathsf{M}}_x + \kappa\hat\sigma \nabla\cdot(\kappa\hat\sigma\nabla\tilde\sigma + \overline{{\boldsymbol{\mathsf{N}}}} )_x.
\end{align}
Taking $L^2$ inner product of \eqref{BB34} with $-\tilde\sigma_x$, we have 
\begin{align}\label{BB35}
-\frac{\mathrm{d}}{\mathrm{d}t} \langle \tilde\sigma_{xt},\tilde\sigma_t\rangle + \|\tilde\sigma_{xt}\|^2 = \int_\Omega \overline{\mathsf{M}}_x \tilde\sigma_x \mathrm{d}\mathbf{x} - \kappa\hat\sigma\int_\Omega \nabla\cdot(\kappa\hat\sigma\nabla\tilde\sigma+\overline{{\boldsymbol{\mathsf{N}}}})_x \tilde\sigma_x \mathrm{d}\mathbf{x}.
\end{align}
For the first integral on the right-hand side of \eqref{BB35}, a direct calculation shows that 
\begin{align}\label{BB36}
\int_\Omega \overline{\mathsf{M}}_x \tilde\sigma_x \mathrm{d}\mathbf{x} = &\,\int_\Omega \tilde\sigma_x \bu \cdot\nabla \tilde\sigma_{xt} \mathrm{d}\mathbf{x} + \kappa \int_\Omega \tilde\sigma_x \tilde\sigma \nabla\cdot\bu_{xt}\mathrm{d}\mathbf{x} + \int_\Omega (\bu_x\cdot\nabla\tilde\sigma_t + \kappa \tilde\sigma_x \nabla\cdot\bu_t) \tilde\sigma_x\mathrm{d}\mathbf{x} \notag \\
&\,+ \int_\Omega (\kappa\tilde\sigma_t\nabla\cdot\bu + \bu_t \cdot\nabla\tilde\sigma)_x \tilde\sigma_x \mathrm{d}\mathbf{x}.
\end{align}
For the first integral on the right-hand side of \eqref{BB36}, using the boundary condition of $\bu$, we have 
\begin{align}\label{BB37}
\int_\Omega \tilde\sigma_x \bu \cdot\nabla \tilde\sigma_{xt} \mathrm{d}\mathbf{x} = - \int_\Omega \tilde\sigma_{xt} \nabla\cdot(\tilde\sigma_x \bu) \mathrm{d}\mathbf{x}.
\end{align}
Since temporal derivative does not change spatial boundary condition, we obtain from \eqref{BB25}--\eqref{BB26}: 
\begin{align*}
&v_{xt}=0,\quad \tilde\sigma_y=0,\quad y=0,1,\ 0\leqslant x \leqslant 1,\\
&u_{yt}=0,\quad \tilde\sigma_x=0,\quad x=0,1,\ 0\leqslant y \leqslant 1,
\end{align*}
which imply 
\begin{equation}\label{BB38}
\tilde\sigma_x\bu_{xt}\cdot\mathbf{n} =\left\{\begin{aligned} 
&\pm \tilde\sigma_x v_{xt} =0,\quad y=0,1,\ 0\leqslant x \leqslant 1,\\
&\pm \tilde\sigma_x u_{xt} =0,\quad x=0,1,\ 0\leqslant y \leqslant 1.
\end{aligned}
\right.
\end{equation}
By \eqref{BB38}, we rewrite the second integral on the right-hand side of \eqref{BB36} as 
\begin{align}\label{BB39}
\kappa\int_\Omega \tilde\sigma_x \tilde\sigma \nabla\cdot\bu_{xt}\mathrm{d}\mathbf{x} = -\kappa\int_\Omega \bu_{xt} \cdot\nabla(\tilde\sigma_x \tilde\sigma) \mathrm{d}\mathbf{x}.
\end{align}
Substituting \eqref{BB37} and \eqref{BB39} into \eqref{BB36}, we obtain 
\begin{align}\label{BB40}
\int_\Omega \overline{\mathsf{M}}_x \tilde\sigma_x \mathrm{d}\mathbf{x} = &\,- \int_\Omega \tilde\sigma_{xt} \nabla\cdot(\tilde\sigma_x \bu) \mathrm{d}\mathbf{x} -\kappa\int_\Omega \bu_{xt} \cdot\nabla(\tilde\sigma_x \tilde\sigma) \mathrm{d}\mathbf{x} + \int_\Omega (\bu_x\cdot\nabla\tilde\sigma_t + \kappa \tilde\sigma_x \nabla\cdot\bu_t) \tilde\sigma_x\mathrm{d}\mathbf{x} \notag \\
&\,+ \int_\Omega (\kappa\tilde\sigma_t\nabla\cdot\bu + \bu_t \cdot\nabla\tilde\sigma)_x \tilde\sigma_x \mathrm{d}\mathbf{x}.
\end{align}
For the second integral on the right-hand side of \eqref{BB35}, by \eqref{N1b} and \eqref{BB38}, we have
\begin{align}\label{BB41}
- \kappa\hat\sigma \int_\Omega \nabla\cdot(\kappa\hat\sigma\nabla\tilde\sigma + \overline{{\boldsymbol{\mathsf{N}}}})_x \tilde\sigma_x \mathrm{d}\mathbf{x} =  - \kappa \hat\sigma\int_{\Omega} \bu_{xt} \cdot \nabla \tilde\sigma_x\mathrm{d}\mathbf{x}.
\end{align}
Substituting \eqref{BB40} and \eqref{BB41} into \eqref{BB35} gives
\begin{align}\label{BB42}
&\,-\frac{\mathrm{d}}{\mathrm{d}t} \langle \tilde\sigma_{xt},\tilde\sigma_t\rangle + \|\tilde\sigma_{xt}\|^2 \notag \\
= &\,- \int_\Omega \tilde\sigma_{xt} \nabla\cdot(\tilde\sigma_x \bu) \mathrm{d}\mathbf{x} -\kappa\int_\Omega \bu_{xt} \cdot\nabla(\tilde\sigma_x \tilde\sigma) \mathrm{d}\mathbf{x} + \int_\Omega (\bu_x\cdot\nabla\tilde\sigma_t + \kappa \tilde\sigma_x \nabla\cdot\bu_t) \tilde\sigma_x\mathrm{d}\mathbf{x} \notag \\
&\,+ \int_\Omega (\kappa\tilde\sigma_t\nabla\cdot\bu + \bu_t \cdot\nabla\tilde\sigma)_x \tilde\sigma_x \mathrm{d}\mathbf{x} - \kappa \hat\sigma \int_{\Omega} \bu_{xt} \cdot \nabla \tilde\sigma_x\mathrm{d}\mathbf{x}.
\end{align}
At this point, we have illustrated how to reallocate the spatial derivatives such that when taking another temporal derivative of \eqref{BB34} and performing spatial integration, no term involving the fourth order derivatives of the unknown functions will appear. To proceed, by applying Sobolev embedding and Gagliardo-Nirenberg inequality, we can derive from \eqref{BB42} that 
\begin{align}\label{BB43}
-\frac{\mathrm{d}}{\mathrm{d}t} \langle \tilde\sigma_{xt},\tilde\sigma_t\rangle + \|\tilde\sigma_{xt}\|^2 \leqslant  C\mathfrak{E}^\frac32 + \frac{\kappa\hat\sigma}{2}\|(\bu_{xt},\nabla \tilde\sigma_x)\|^2.
\end{align}
In a similar fashion, it can be shown that 
\begin{align}\label{BB44}
-\frac{\mathrm{d}}{\mathrm{d}t} \langle \tilde\sigma_{yt},\tilde\sigma_t\rangle + \|\tilde\sigma_{yt}\|^2 \leqslant C\mathfrak{E}^\frac32 + \frac{\kappa\hat\sigma}{2}\|(\bu_{yt},\nabla \tilde\sigma_y)\|^2.
\end{align}
Combining \eqref{BB34}, \eqref{BB42}, and \eqref{BB43}, we obtain 
\begin{align}\label{BB45}
-\frac{\mathrm{d}}{\mathrm{d}t} \langle\tilde\sigma_t,\tilde\sigma\rangle_{1} + \|\tilde\sigma_{t}\|_{H^1}^2 \leqslant C\mathfrak{E}^\frac32 + \frac{\kappa\hat\sigma}{2} \|(\bu_{t},\nabla \tilde\sigma)\|_{H^1}^2,
\end{align}
where $\langle\cdot,\cdot\rangle_{1}$ denotes the standard inner product in $H^1$. By taking temporal derivative of \eqref{BB34} and repeating the above arguments, we have the following:
\begin{align}\label{BB46}
-\frac{\mathrm{d}}{\mathrm{d}t} \langle\tilde\sigma_{tt},\tilde\sigma_t\rangle_{1} + \|\tilde\sigma_{tt}\|_{H^1}^2 \leqslant C\mathfrak{E}^\frac32 + \frac{\kappa\hat\sigma}{2} \|(\bu_{tt},\nabla \tilde\sigma_t)\|_{H^1}^2.
\end{align}
The combination of \eqref{BB45} and \eqref{BB46} yields 
\begin{align}\label{BB47}
-\frac{\mathrm{d}}{\mathrm{d}t} \Big(\sum_{k=1}^2\langle \partial_t^k\tilde\sigma,\partial_t^{k-1}\tilde\sigma\rangle_{1}\Big) + \sum_{k=1}^2\|\partial_t^k\tilde\sigma\|_{H^1}^2 \leqslant C\mathfrak{E}^\frac32 + \frac{\kappa\hat\sigma}{2} \|(\bu_{t},\bu_{tt},\nabla \tilde\sigma,\nabla \tilde\sigma_t)\|_{H^1}^2.
\end{align}
Note that according to \eqref{BB7} and \eqref{BB9}, it holds that 
\begin{align}\label{BB48}
\|(\nabla\tilde\sigma,\nabla\tilde\sigma_t)\|_{H^1}^2 \leqslant C\big(\|(\bu_{tt},\bu_t,\bu)\|_{H^1}^2+\|(\nabla\tilde\phi,\nabla\tilde\phi_t)\|_{H^1}^2+\mathfrak{E}^2\big).
\end{align}
By \eqref{BB48}, we update \eqref{BB47} as
\begin{align}\label{BB49}
-\frac{\mathrm{d}}{\mathrm{d}t} \Big(\sum_{k=1}^2\langle \partial_t^k\tilde\sigma,\partial_t^{k-1}\tilde\sigma\rangle_{1}\Big) + \sum_{k=1}^2\|\partial_t^k\tilde\sigma\|_{H^1}^2 \leqslant C\mathfrak{E}^\frac32 + \mathsf{C}_2\Big(\sum_{k=0}^2\|\partial_t^k\bu\|_{H^1}^2+\sum_{k=0}^1\|\partial_t^k\nabla\tilde\phi\|_{H^1}^2\Big).
\end{align}
This completes the recovery of density dissipation. Next, we deal with the estimate of $\tilde\phi$.

\vspace{.1 in}

{\bf Estimation of $\tilde\phi$}: Using \eqref{N1c} and Poincar\'e's inequality, we can show that 
\begin{align}\label{BB50}
\tau \frac{\mathrm{d}}{\mathrm{d}t}\Big(\sum_{k=0}^2\|\partial_t^k\tilde\phi\|_{H^1}^2\Big) + d \sum_{k=0}^2\|(\partial_t^k\nabla\tilde\phi,\partial_t^k\Delta\tilde\phi)\|^2 \leqslant Cd^{-1}\mathfrak{E},
\end{align}
where we dropped the non-negative terms involving $a$. Since the calculations are straightforward, we omitted the details for brevity. Applying \eqref{AA14} to $\partial_t^k\nabla\tilde\phi_t$, we have 
$$
\|\partial_t^k\nabla\tilde\phi\|_{H^1}^2 \leqslant C\big( \|\partial_t^k\Delta\tilde\phi\|^2 + \|\partial_t^k\nabla\tilde\phi\|^2\big),\quad k=0,1,2.
$$ 
Then we  update \eqref{BB50} as 
\begin{align}\label{BB51}
\tau \frac{\mathrm{d}}{\mathrm{d}t}\Big(\sum_{k=0}^2\|\partial_t^k\tilde\phi\|_{H^1}^2\Big)+\mathsf{C}_3d\sum_{k=0}^2\|\partial_t^k\nabla\tilde\phi\|_{H^1}^2  \leqslant Cd^{-1}\mathfrak{E}.
\end{align}
Now we are ready to couple different energy estimates to yield the desired results.

\vspace{.1 in}

{\bf Closing of energy estimates}: Let $\mathsf{L}_2$ be a constant such that $\mathsf{L}_2>4\max\{1,\alpha^{-1}\mathsf{C}_2\}$. Then the operation $\mathsf{L}_2\times\eqref{BB31}$+\eqref{BB49} yields
\begin{align}\label{BB52}
\frac{\mathrm{d}}{\mathrm{d}t}\mathtt{X}_2 + \mathtt{Y}_2 \leqslant  C\mathfrak{E}^\frac32 + C\sum_{k=0}^2\|\partial_t^k\nabla\tilde\phi\|_{H^1}^2,
\end{align}
where the quantities on the left-hand side are defined by 
\begin{align*}
\mathtt{X}_2&\triangleq \mathsf{L}_2\sum_{k=0}^2\|(\partial_t^k\tilde\sigma,\partial_t^k\bu)\|_{H^1}^2 -\sum_{k=1}^2\langle \partial_t^k\tilde\sigma,\partial_t^{k-1}\tilde\sigma\rangle_{1},\\
\mathtt{Y}_2&\triangleq (\alpha\mathsf{L}_2-\mathsf{C}_2)\sum_{k=0}^2\|\partial_t^k\bu\|_{H^1}^2 + \sum_{k=1}^2\|\partial_t^k\tilde\sigma\|_{H^1}^2.
\end{align*}
By the choice of $\mathsf{L}_2$, it is easy to see that 
\begin{align}\label{BB53}
\mathtt{X}_2&\cong \sum_{k=0}^2\|(\partial_t^k\tilde\sigma,\partial_t^k\bu)\|_{H^1}^2.
\end{align}
In addition, because of Poincar\'e's inequality and \eqref{N1b}, we have
$$
\|\tilde\sigma\|^2 \leqslant C\|\nabla\tilde\sigma\|^2 \leqslant C\big(\|(\bu_t,\bu)\|^2+\mathfrak{E}^2+\|\nabla\tilde\phi\|^2\big),
$$
which implies
\begin{align}\label{BB54}
\mathtt{Y}_2&\gtrsim \sum_{k=0}^2\|(\partial_t^k\tilde\sigma,\partial_t^k\bu)\|_{H^1}^2 - \mathfrak{E}^2 - \|\nabla\tilde\phi\|^2.
\end{align}
By \eqref{BB53} and \eqref{BB54}, we upgrade \eqref{BB52} as 
\begin{align}\label{BB55}
\frac{\mathrm{d}}{\mathrm{d}t}\mathtt{X}_2 + C\mathtt{X}_2 \leqslant  C\mathfrak{E}^\frac32 + C\sum_{k=0}^2\|\partial_t^k\nabla\tilde\phi\|_{H^1}^2.
\end{align}
Similar to \eqref{AA52}, we have 
\begin{align}\label{BB56}
\frac{\mathrm{d}}{\mathrm{d}t} \mathtt{V}_2 + \alpha \mathtt{V}_2 \leqslant C\mathfrak{E}^{\frac32}.
\end{align}
Combining \eqref{BB55} and \eqref{BB56}, we get
\begin{align}\label{BB57}
\frac{\mathrm{d}}{\mathrm{d}t} (\mathtt{X}_2+\mathtt{V}_2) + C(\mathtt{X}_2+\mathtt{V}_2) \leqslant  C\mathfrak{E}^\frac32 + \mathsf{C}_4\sum_{k=0}^2\|\partial_t^k\nabla\tilde\phi\|_{H^1}^2.
\end{align}
When $d$ is sufficiently large such that $\mathsf{C}_3 d > \mathsf{C}_4$, by adding \eqref{BB51} and \eqref{BB57}, we obtain
\begin{align}\label{BB58}
\frac{\mathrm{d}}{\mathrm{d}t} (\mathtt{X}_2+\mathtt{V}_2+\mathtt{Z}) + C(\mathtt{X}_2+\mathtt{V}_2+\overline{\mathtt{Z}}) \leqslant  C\big(\mathfrak{E}^\frac12 +d^{-1}\big)\mathfrak{E},
\end{align}
where the new quantities on the left-hand side are defined by 
\begin{align*}
\mathtt{Z}(t) \triangleq \tau \sum_{k=0}^2\|\partial_t^k\tilde\phi\|_{H^1}^2 ,\quad \overline{\mathtt{Z}}(t) \triangleq (\mathsf{C}_3d-\mathsf{C}_4)\sum_{k=0}^2\|\partial_t^k\nabla\tilde\phi\|_{H^1}^2.
\end{align*}
In view of \eqref{BB5} and Poincar\'e's inequality, we see that 
\begin{align}\label{BB58a}
\mathtt{W}(t)\lesssim \mathtt{Z}(t)\lesssim \overline{\mathtt{Z}}(t).
\end{align}
Moreover, it is apparent that $\mathtt{U}_2(t) \lesssim \mathtt{X}_2(t)$. Then, by \eqref{BB21}, we update \eqref{BB58} as 
\begin{align}\label{BB59}
\frac{\mathrm{d}}{\mathrm{d}t} (\mathtt{X}_2 + \mathtt{V}_2 + \mathtt{Z}) + C(\mathtt{X}_2+\mathtt{V}_2+\overline{\mathtt{Z}}) \leqslant  C\big(\mathfrak{E}^\frac12 + d^{-1}\big)(\mathtt{X}_2+\mathtt{V}_2+\overline{\mathtt{Z}}).
\end{align}
When $\mathfrak{E}$ is sufficiently small and $d$ is sufficiently large, we obtain from \eqref{BB59}:
 \begin{align}\label{BB60}
\frac{\mathrm{d}}{\mathrm{d}t} (\mathtt{X}_2+\mathtt{V}_2+\mathtt{Z}) + C(\mathtt{X}_2+\mathtt{V}_2+\overline{\mathtt{Z}}) \leqslant  0.
\end{align}
Integrating \eqref{BB60} with respect to $t$, we get the following estimate:
\begin{align}\label{BB61}
\mathcal{E}(t) + \int_0^t \mathfrak{D}(\tau)\mathrm{d}\tau \leqslant C\mathcal{E}(0),
\end{align}
where the quantities are defined by
\begin{align*}
\mathcal{E}\triangleq \sum_{k=0}^2\|(\partial_t^k\tilde\sigma,\partial_t^k\bu,\partial_t^k\tilde\phi)\|_{H^1}^2 + \mathtt{V}_2,\qquad 
\mathfrak{D}\triangleq \sum_{k=0}^2\|(\partial_t^k\tilde\sigma,\partial_t^k\bu,\partial_t^k\nabla\tilde\phi)\|_{H^1}^2 + \mathtt{V}_2.
\end{align*}
Since $\mathcal{E} \geqslant \mathtt{U}_2 + \mathtt{V}_2 + \mathtt{W}$, by 
\eqref{BB21} and \eqref{BB61}, we deduce that
\begin{align*}
\mathfrak{E}(t) \leqslant C\mathcal{E}(t) \leqslant C\mathcal{E}(0).
\end{align*}
By virtue of \eqref{BB58a} and \eqref{BB60}, we know that $\mathcal{E}(t)$ decays exponentially in time. Hence, $\mathfrak{E}(t)$ decays exponentially in time. Moreover, since $\mathcal{E}\lesssim \mathfrak{D}$, we have 
\begin{align*}
\int_0^t \mathfrak{E}(\tau)\mathrm{d}\tau \leqslant C\mathcal{E}(0).
\end{align*}
Furthermore, since by Poincar\'e's inequality, $\|\tilde\phi_{tt}\|_{H^2} \lesssim \|\nabla\tilde\phi_{tt}\|_{H^1}$, we conclude that 
\begin{align*}
\int_0^t \|\tilde\phi_{tt}(\tau)\|_{H^2}^2\mathrm{d}\tau \leqslant C\mathcal{E}(0).
\end{align*}
This completes the proof of Theorem \ref{main1}  for $\gamma > 1$.

\subsection{Proof of Theorem \ref{main1} for $\gamma=1$} The proof follows similar arguments to the $\gamma>1$ case, with modifications required for the symmetrization of the fluid component. We begin by multiplying \eqref{HP1} by $A_0$ and \eqref{HP2} by $\rho$, yielding the symmetric system:
\begin{subequations}\label{Z1}
\begin{alignat}{2}
A_0\rho_t + A_0\bu\cdot\nabla\rho + A_0\rho\nabla\cdot\bu&=0, \label{Z1a}\\
\rho^2\bu_t+\rho^2\bu\cdot\nabla\bu+A_0\rho\nabla\rho&=-\alpha\rho^2\bu+\beta\rho^2\nabla\phi. \label{Z1b}
\end{alignat}
\end{subequations}
Introducing the perturbation variables $\tilde\rho=\rho-\hat\rho_0$ and $\tilde\phi=\phi-\hat\phi(t)$ (with $\hat\rho_0$ and $\hat\phi(t)$ defined as before), we rewrite \eqref{Z1} together with \eqref{HP3} as 
\begin{subequations}\label{Z2}
\begin{alignat}{3}
A_0\tilde\rho_t + A_0\bu\cdot\nabla\tilde\rho + A_0\tilde\rho\nabla\cdot\bu+A_0\hat\rho_0\nabla\cdot\bu&=0, \label{Z2a}\\
(\tilde\rho+\hat\rho_0)^2\bu_t+(\tilde\rho+\hat\rho_0)^2\bu\cdot\nabla\bu+A_0\tilde\rho\nabla\tilde\rho+A_0\hat\rho_0\nabla\tilde\rho&=-\alpha(\tilde\rho+\hat\rho_0)^2\bu+\beta(\tilde\rho+\hat\rho_0)^2\nabla\tilde\phi, \label{Z2b}\\
\tau\tilde\phi_t-d\Delta \tilde\phi + a \tilde\phi &=b\tilde\rho. \label{Z2c}
\end{alignat}
\end{subequations}
The energy reduction and estimation of reduced energy require analogous calculations as in the $\gamma>1$ case. For density dissipation recovery, we derive the governing equation:
\begin{align}
\tilde\rho_{tt}+\nabla\cdot(\tilde\rho\bu)_t-A_0\hat\rho_0\nabla\cdot\big[(\tilde\rho+\hat\rho_0)^{-1}\nabla\tilde\rho\big] = \hat\rho_0\nabla\cdot(\bu\cdot\nabla\bu+\alpha\bu-\beta\nabla\tilde\phi),
\end{align}
obtained by combining \eqref{Z2a} and \eqref{Z2b}. The remaining analysis, including the estimation of $\tilde\phi$ and $\nabla\times \bu$, as well as final energy estimate closure, follows, with the necessary modifications, from the $\gamma>1$ case. We omit these technical details for brevity. This completes the proof of Theorem \ref{main1} for $\gamma=1$, and the paper ends here.

\section{Conclusion and Outlook}

This paper presents a rigorous mathematical analysis of a PDE system modeling vasculogenesis in a two-dimensional bounded domain. The primary contributions are twofold:
\begin{enumerate}
\item[1.] For the system with quadratic pressure and Dirichlet boundary condition for the chemical attractant, explicit formulas for non-constant steady-state solutions are derived on the unit square using eigenfunction expansion. Their existence is guaranteed under a condition requiring a sufficiently large diffusion coefficient.

\item[2.] The asymptotic stability of the constructed steady states is proven under two scenarios:
\begin{itemize}
\item For the fast-relaxation limit of the quadratic pressure system with Dirichlet condition for $\phi$, provided the initial perturbation is small and the diffusion coefficient is sufficiently large, the solution converges to the steady state exponentially in time.

\item For the general $\gamma$-law pressure with Neumann condition for $\phi$, the constant steady state determined by the initial and boundary conditions is shown to be locally exponentially stable under similar assumptions on the diffusion coefficient  and initial data.
\end{itemize}
\end{enumerate}
The analysis overcomes significant technical challenges, particularly the coupling between hyperbolic and parabolic/elliptic components, and the handling of boundary terms in the absence of direct boundary conditions for the density. Key techniques include sound-speed transformation, reduction of spatial derivatives to temporal ones, energy estimates leveraging Poincar\'e inequalities and elliptic regularity, and careful treatment of domain geometry and boundary conditions.

Several promising directions emerge from this study:
\begin{enumerate}
\item[1.] The most immediate challenge is to extend the analysis to three-dimensional domains, which are biologically most relevant. This requires validating fundamental analytical tools on 3D polyhedral domains and potentially developing new techniques to handle the increased complexity of the system's structure in higher dimensions.

\item[2.] The stability results for the Dirichlet case are currently limited to the fast-relaxation limit. Establishing stability for the full system remains an open and mathematically demanding problem due to the boundary condition.

\item[3.] The requirement of a sufficiently large diffusion coefficient is crucial for the current analysis but may be a technical rather than biological necessity. Future work could aim to establish stability under milder or more physically realistic assumptions on the parameters.

\item[4.] While the square domain facilitates the use of eigenfunction expansion, extending the results to general smooth or polygonal domains would enhance the model's applicability. Investigating other biologically plausible boundary conditions is also of interest.

\item[5.] Numerical studies comparing the theoretical predictions with simulations of the full time-dependent PDE system could validate the analysis and explore regimes beyond the theoretical assumptions.

\item[6.] A future interdisciplinary goal is to use the analytically tractable steady states and stability criteria as a framework for interpreting {\it in vitro} experiments or for designing control strategies in tissue engineering applications related to vascular network formation.
\end{enumerate}

\section*{Acknowledgements} 

Support of this work came in part from the Fundamental Research Funds for Central Universities of China No.\,3072024WD2401 (K. Zhao).


\begin{thebibliography}{99}
\bibitem{1} F. Ambrosi, D. Bussolino, and L. Preziosi, A review of vasculogenesis models, {\it J. Theo. Med.}, {6}(1): 1--19, 2005.

\bibitem{Bar} C. Bardos, F. Di Plinio, and R. Temam, The Euler equations in planar nonsmooth convex domains, {\it J. Math. Anal. Appl.}, {407}: 69--89, 2013.

\bibitem{2} F. Berthelin, D. Chiron, and M. Ribot, Stationary solutions with vacuum for a one-dimensional chemotaxis model with nonlinear pressure, {\it  Comm. Math. Sci.}, {14}(1): 147--186, 2016.

\bibitem{3} P. Carmeliet, Mechanisms of angiogenesis and arteriogenesis, {\it Nature Med.}, {6}(4): 389--395, 2000.

\bibitem{4} J. Carrillo, X. Chen, Q. Wang, Z. Wang, and L. Zhang, Phase transitions and bump solutions of the Keller-Segel model with volume exclusion, {\it SIAM J. Appl. Math.}, {80}: 232--261, 2020.

\bibitem{5} P. Chavanis and C. Sire, Kinetic and hydrodynamic models of chemotactic aggregation, {\it Physica A}, {384}: 199--222, 2007.

\bibitem{6} T. Crin-Barat, Q. He, and L. Shou, The hyperbolic-parabolic chemotaxis system for vasculogenesis: Global dynamics and relaxation limit toward a Keller-Segel model, {\it SIAM J. Math. Anal.}, {55}(5): 4445--4492, 2023.

\bibitem{7}  M. Di Francesco and D. Donatelli, Singular convergence of nonlinear hyperbolic chemotaxis systems to Keller-Segel type models, {\it Disc. Cont. Dyna. Syst.-Ser. B}, {13}: 79--100, 2010.

\bibitem{8} A. Gamba, D. Ambrosi, A. Coniglio, A. de Candia, S. Di Talia, E. Giraudo, G. Serini, L. Preziosi, and F. Bussolino, Percolation, morphogenesis, and Burgers dynamics in blood vessels formation, {\it Phys. Rev. Let.}, {90}(11): 118101--118104, 2003.

\bibitem{Gri} P. Grisvard, {\it Elliptic Problems in Nonsmooth Domains}, Monographs and Studies in Mathematics, { 24}, Pitman (Advanced Publishing Program), Boston, MA, 1985.

\bibitem{9}  G. Helmlinger, M. Endo, N. Ferrara, L. Hlatky, and R. Jain, Formation of endothelial cell networks, {\it Nature}, {405}: 139--141, 2000.

\bibitem{HHYZ} S. Hertrich, T. Huang, D. Y\'epez, and K. Zhao, Stationary solutions with vacuum for a hyperbolic-parabolic chemotaxis model in dimension two, {\it Nonlinear Anal. RWA}, 88: 104489, 2026. 

\bibitem{10}  G. Hong, H. Peng, Z. Wang, and C. Zhu, Nonlinear stability of phase transition steady states to a hyperbolic–parabolic system modeling vascular networks, {\it J. London Math. Soc.}, {103}: 1480--1514, 2021.

\bibitem{12}  R. Kowalczyk, A. Gamba, and L. Preziosi, On the stability of homogeneous solutions to some aggregation models, {\it Disc. Cont. Dyna. Syst.-Ser. B}, {4}: 203--220, 2004.

\bibitem{Lac} C. Lacave, E. Miot, and C. Wang, Uniqueness for the 2-D Euler equations on domains with corners, {\it Indiana Univ. Math. J.}, {63}: 1725--1756, 2014.

\bibitem{14}  Q. Liu, H. Peng, and Z. Wang, Asymptotic stability of diffusion waves of a quasi-linear hyperbolic-parabolic model for vasculogenesis, {\it SIAM J. Math. Anal.}, {54}: 1313--1346, 2022.

\bibitem{15}  Q. Liu, H. Peng, and Z. Wang, Convergence to nonlinear diffusion waves for a hyperbolic-parabolic chemotaxis system modeling vasculogenesis, {\it J. Diff. Equ.}, {314}: 251--286, 2022.

\bibitem{MN} A. Matsumura and T. Nishida, Initial boundary value problems for the equations of motion of compressible viscous and heat-conductive fluids, {\it Comm. Math. Phys.}, {89}: 445--464, 1983.

\bibitem{16}  R. Natalini, M. Ribot, and M. Twarogowska, A numerical comparison between degenerate parabolic and quasilinear hyperbolic models of cell movements under chemotaxis, {\it J. Sci. Comp.}, {63}: 654--677, 2015.

\bibitem{17}  R. Pan and K. Zhao, The 3D compressible Euler equations with damping in a bounded domain, {\it J. Diff. Equ.}, {246}: 581--596, 2009.

\bibitem{18}  H. Peng and K. Zhao, On a hyperbolic-parabolic chemotaxis system, {\it Math. Biosci. Engin.}, {20}(5): 7802--7827, 2023.

\bibitem{19}  C. Di Russo, Analysis and numerical approximations of hydrodynamical models of biological movements, {\it Rend. Mat. Appl.}, {32}: 117--367, 2012.

\bibitem{20}  C. Di Russo and A. Sepe, Existence and asymptotic behavior of solutions to a quasi-linear hyperbolic-parabolic model of vasculogenesis, {\it SIAM J. Math. Anal.}, {45}: 748--776, 2013.

\bibitem{21}  S. Schochet, The compressible Euler equations in a bounded domain: Existence of solutions and the incompressible limit, {\it Comm. Math. Phys.}, {104}: 49--75, 1986.

\bibitem{22} Y. Xiao and Z.-P. Xin, On the vanishing viscosity limit for the 3D Navier-Stokes equations with a slip boundary condition, {\it Comm. Pure Appl. Math.}, {60}(7): 1027--1055, 2007.
\end{thebibliography}
\end{document}